\documentclass[a4paper,11pt]{article}
\usepackage{nicefrac}
\usepackage{todonotes}
\definecolor{refkey}{rgb}{.35,.75,0}
\definecolor{labelkey}{rgb}{.15,.55,0}

\usepackage{mathtools}
\usepackage[utf8]{inputenc}
\usepackage{pgfplots}\pgfplotsset{compat=1.18}
\usepackage{tikz}

\usepackage{float}

\usepackage{epsfig}
\usepackage{amssymb,amsmath,amsfonts}

\catcode`@=11
\@addtoreset{equation}{section}
\catcode`@=12
\newcommand{\dx}{\hbox{d}x}
\newcommand{\dy}{\hbox{d}y}
\newcommand{\ddt}{\mathrm{\Delta t}}
\newcommand{\ddx}{\mathrm{\Delta x}}
\newcommand{\ddy}{\mathrm{\Delta y}}

\newcommand{\R}{{\mathbb R}}

\newcommand{\Z}{\mathbb Z}

\newcommand{\Rd}{\R^d}

\newcommand{\T}{{\cal{T}}}

\newcommand{\bee}{\begin{equation}}
\newcommand{\ene}{\end{equation}}
\newcommand{\bes}{\begin{section}}
\newcommand{\ens}{\end{section}}

\def\refp#1{(\ref{#1})}

\newcommand{\OSC}{\mathrm{I}}
\newcommand{\ENO}{\ensuremath{\mathsf{ENO}}}
\newcommand{\WENO}{\ensuremath{\mathsf{WENO}}}
\newcommand{\CWENO}{\ensuremath{\mathsf{CWENO}}}
\newcommand{\CWENOZ}{\ensuremath{\mathsf{CWENOZ}}}

\newcommand{\Ogrande}{\mathcal{O}}
\newcommand{\stencil}{\mathcal{S}}

\newtheorem{Theorem}{Theorem}[section]
\newtheorem{Definition}[Theorem]{Definition}

\begin{document}
\title{A CWENO large time-step scheme for Hamilton--Jacobi equations}
\author{E. Carlini, R. Ferretti, S. Preda, and M. Semplice}
\date{}
\maketitle
\begin{abstract}

\end{abstract}
We propose a high order numerical scheme for time-dependent first order Hamilton--Jacobi--Bellman equations. In particular we propose to combine a semi-Lagrangian scheme with a  Central Weighted Non-Oscillatory reconstruction. We prove a convergence result in the case of state- and time-independent Hamiltonians.

Numerical simulations are presented in space dimensions one and two, also for more general state- and time-dependent Hamiltonians, demonstrating superior performance in terms of CPU time gain compared with a semi-Lagrangian scheme coupled with Weighted Non-Oscillatory reconstructions.

%\vskip0.5cm\noindent
%\vskip0.5cm\noindent
%\vskip1cm\noindent

\section{Introduction}
\label{sec:setting}

In this paper, we address the numerical approximation of the following hyperbolic Hamilton--Jacobi--Bellman (HJB) equation:
\begin{equation}\label{eq:HJ}
\begin{cases}
v_t(t,x) + H(t,x,Dv(t,x)) = 0, & \text{for } t,x \in (0,T) \times \mathbb{R}^d, \\
v(0,x) = v_0(x), & \text{for } x \in \mathbb{R}^d,
\end{cases}
\end{equation}
where $v: (0,T) \times \mathbb{R}^d \to \mathbb{R}$, $Dv$ stands for the spatial gradient, and $H: (0,T) \times \mathbb{R}^d \times \mathbb{R}^d \to \mathbb{R}$.

Numerous schemes have been proposed to approximate \eqref{eq:HJ}, but only a small number are aimed at high order accuracy.
%% Finite Difference Schemes
High order finite difference schemes based on $\ENO$ (Essentially Non-Oscillatory) reconstruction, defined in \cite{HEOC87}, were proposed in \cite{OS91} and extended to unstructured grids in \cite{A94}. Second order Godunov-type schemes based on global projection operators are discussed in \cite{LT00}. Weighted $\ENO$ ($\WENO$) schemes were introduced in \cite{JiangShu:96,LOC94} and combined with finite difference schemes for Hamilton--Jacobi (HJ) equations in \cite{JP00}. A fifth order central scheme based on $\WENO$ reconstruction has been developed in \cite{BL03}.

%% Semi-Lagrangian Schemes
Over the past two decades, semi-Lagrangian (SL) schemes have been employed to discretize the HJ equation. In \cite{FF94}, a high order semi-discrete SL scheme is proposed to discretize stationary HJB equations, while in \cite{FR02}, SL schemes are applied to evolutionary HJB equations. A high order SL scheme for HJB equations is presented in \cite{CFR05} by combining the SL technique with $\WENO$ reconstructions.

%% Filtered Schemes
A different approach to achieve high order accuracy involves filtered schemes. In \cite{BFS16,FPT20}, filtered schemes are combined with monotone finite difference schemes to discretize first order evolutionary HJ equations.

For an overview of numerical methods for first order Hamilton--Jacobi equations, we also refer to the book \cite{FF14}.

We will develop our theory for the HJB equation related to a finite horizon optimal control. More precisely, given a compact \(A \in \mathbb{R}^m\), 
%for \(m \leq N\), 
%\todo{Cos'è $N$? Togliamo $m\leq N$?}
let the set \(\mathcal{A} = \{\alpha:[0,T] \to A, \text{ measurable}\}\) denote the admissible controls. Given the running cost \(f_C: (0,T) \times \mathbb{R}^d \times A \to \mathbb{R}\) and the dynamics \(f_D: (0,T) \times \mathbb{R}^d \times A \to \mathbb{R}^d\) for the system
\begin{equation}\label{eq:ode}
\dot{y}(s) = f_D(t-s, y(s), \alpha(s)), \quad s \in (0,t], \quad y(0) = x,
\end{equation}
we consider the following Hamiltonian
\begin{equation}\label{eq:Hamiltonian}
H(t,x,p) = \max_{a \in A}(-f_D(t,x,a) \cdot p - f_C(t,x,a)).
\end{equation}

Let us assume that
\begin{itemize}
%\item[\bf(H1)] $A$ is a compact set of $\mathbb{R}^m$, for $m\leq N$
\item[\bf(H1)] $f_C,f_D$ are continuous and bounded. Moreover, for every $a\in A$, the functions $f_C(\cdot,\cdot,a)$, $f_D(\cdot,\cdot,a)$ are Lipschitz continous, with Lipschitz constants
independent of $a \in A$.
\item[\bf(H2)] $v_0$ is Lipschitz continous and bounded.
\end{itemize}
Under Assumptions {\bf(H1)--\bf(H2)}, problem \eqref{eq:HJ} admits a unique viscosity solution $v$, which is Lipschitz continuous and bounded. Moreover, the Dynamic Programming Principle holds (see \cite[Chapter 3]{BCD97}, i.e., for any $h>0$
\begin{equation}\label{eq:DP}
v(t,x)= \underset {\alpha \in \mathcal A}\inf\left\{ \int_0^h f_C(t-s,y_{x,t}(s),\alpha(s)) {\rm{d}}s + v(y_{x,t}(h),t-h)\right\},
\end{equation}
where $y_{x,t}$ denotes the solution to \eqref{eq:ode}.

%With this setting, the solution to \eqref{eq:HJ} is the value function of a finite
%horizon optimal control problem and the following representation formula holds
%\begin{equation}\label{eq:RP}
%v(x,t)=\underset{\alpha \in \A}{\inf }\left\{ \int_0^t f_C(s,y(s),\alpha(s))\mathrm{d} s+v_0(y(t)) \right\}
%\end{equation}
%where $y(t)$ denotes the solution to \eqref{eq:ode}.
%{\color{red} inserire HP sui dati e richianare risulati esistenza e regolarita\\
%inserire formula programmazione dinamica per introdurre schema 2.1 ELISABETTA\\}

A numerical scheme based on the previous description involves a minimization procedure over the set of controls; in turn, each function evaluation requires to approximate the cost integral and to interpolate  the solution at the previous time step at the foot $y_{x,t}(h)$ of the characteristics. In this respect, $\WENO$ reconstructions, especially in higher space dimensions, are not efficient, due to fact that the linear weights of $\WENO$ depend on the reconstruction point.

A novel paradigm for non-oscillatory reconstruction operators has been introduced in \cite{LPR:00:SIAMJSciComp}, where the authors suggested to blend, in a $\WENO$-like fashion, polynomials of different degrees, allowing to overcome some difficulties of non-existence, non-positivity and dependence on the reconstruction point of the $\WENO$ linear weights. The idea has been further developed into the so-called $\CWENO$ reconstruction and exploited in more spatial dimensions, also in the case of adaptive mesh refinement and non-uniform grids 
\cite{BGFB:2020,ZhuQiu:2016,Baeza:19:CWENOglobalaverageweight,ZhouCai:08,ADER_CWENO}.
The technique has also been exploited in finite difference schemes for Hamilton-Jacobi equations on Cartesian meshes via dimensional splitting \cite{ZhuQiu:2017:HJ,ZhengShuQiu:2019:HJ} and on general meshes \cite{ZhuQiu:2020:triHJ}.
General results for establishing the convergence order of a $\CWENO$ reconstruction have been presented in \cite{CPSV:cweno,CSV19:cwenoz,SV:2020:CWAO}.

One of the main advantages of $\CWENO$ over the traditional $\WENO$ is that the central approach provides a reconstruction polynomial that is defined everywhere in the reconstruction cell and that can be evaluated later, with no essential extra cost, at many different reconstruction points. This is guaranteed by the independence of the linear weights from the reconstruction point. Furthermore, as shown in \cite{SCR:cwenoAMR,ADER_CWENO,SC:19:wb2d,CPSV:cweno}, the $\CWENO$ approach allows to avoid the dimensional splitting procedure and this is advantageous also on Cartesian meshes when the number of reconstruction points per cell is high.

In this paper we want to exploit the positive features of $\CWENO$ to obtain a high order SL scheme that is more efficient than the one of \cite{CFR05}, which are based on $\WENO$.
The rest of the paper is organized as follows. In \S\ref{sec:numericalScheme} the general principles of the SL scheme are recalled, whereas in \S \ref{sec:cweno} a $\CWENO$ interpolation from point values is constructed in one and two space dimensions. \S\ref{sec:convergence} proves a convergence result in the framework of \cite{F01}, and finally \S\ref{sec:tests} provides an extensive number of numerical examples to validate the scheme.

\section{Numerical scheme}
\label{sec:numericalScheme}
To obtain an approximate version of \eqref{eq:HJ}, we need first to discretize the control problem in time. %, as formulated in \cite{FF94} for high order approximations.
Let $\ddt > 0$ be a time step, $t_n = n\ddt$ a uniform time grid with $n = 0, \ldots, N_T$, where $N_T = \lceil \frac{T}{\ddt} \rceil$, and consider the Dynamic Programming Principle on a single time step $[t_n,t_{n+1}]$, by choosing in \eqref{eq:DP} $h=\ddt$ and $t=t_{n+1}$.

When using a $\nu$-stages scheme for \eqref{eq:ode}, we discretize an admissible control $\alpha\in\mathcal{A}$ via a sequence 
$\underline{a} = (a_1, a_2, \ldots, a_{\nu})\in A^\nu$. 
Then the foot of the characteristic corresponding to the control $\underline{a}$ is
%We use a $\nu$-stages scheme to approximate   the solution to \eqref{eq:ode} for $t=t_{n+1}$ and for a time  step $\ddt$, 
%\[
%y_i^{n}(\underline{a}) = x_i + \ddt \Phi_D(t_n, x_i, \underline{a}),
%\]
\[
y^{n}(x,\underline{a}) \simeq y_{x,t_{n+1}}(\ddt)% x+\int_{t_n}^{t_{n+1}} f_D(t_{n+1}-s,y_{x,t_{n+1}}(s),\alpha(s))\, {\rm d} s,
\]
 for any $n = 0, \ldots, N_T-1$, $x \in \R^d$ and $\underline{a}  \in A^{\nu}$. 
% The increment function $\Phi_{D}(t_n, x_i, \underline{a})$ is such that
%\[
%\lim_{
%\substack{\ddt\to 0 \\ t_n\to t}
%}
%\Phi_{D}(t_n, x, \underline{a}) = f_D(x, t, \underline{a}),
%\]
%for any $x\in \mathbb{R}^d$.

Moreover, a suitable quadrature formula based on $\nu$ quadrature nodes is introduced to discretize the integral term related to the running cost,
 \[
C^{n}(x,\underline{a})\simeq \int_{0}^{\ddt} f_C(t_{n+1}-s,y_{x,t_{n+1}}(s),\alpha(s))\, {\rm d} s,
\] 
 for any $n = 0, \ldots, N_T-1$, $x \in \R^d$ and $\underline{a}  \in A^{\nu}$ .
 
Let us now introduce the space grid with space step $\ddx$ so that $x_i=i\ddx$ with $i \in \Z^d$ and denote $y_i^{n}(\underline{a})=y^{n}(x_i,\underline{a})$ and $C^{n}_i(\underline{a})=C^{n}(x_i,\underline{a})$.
We approximate the solution $v(x,t)$ of \eqref{eq:HJ} by a discrete function $u^n_i\simeq v(x_i,t_n)$ defined on the space-time grid, computed by the following iterative scheme, for $n=0,\dots,N_T-1$
% which is second order in time and use a central Weno reconstruction is space of {\color{red}order 2-3} 
\begin{equation}\label{eq:SL}
\begin{cases}
u^{n+1}_i=\underset{\underline{a} \in A^\nu}\min \{R[u^n](y^{n}_i(\underline{a}))+C^{n}_i(\underline{a})\} & i\in \Z^d,\\
u^{0}_i=v_0(x_i), & i\in \Z^d,
\end{cases}
\end{equation}
where $R[u^n](x)$ denotes the spatial reconstruction at $x\in \Rd$ of the numerical solution $u^n=(u^n_i)_{i\in\Z^d}$.

%where $\underline{a}=(a_1,a_2,...,a_{\nu})$, $a_j\in A$ $\forall j=1,2,...,\nu$, $y^n_i(\underline{a})\simeq y_{x_i,t_{n+1}}(\ddt)$, i.e. $y^n_i(\underline{a})$ it is the numerical solution of the ODE \eqref{eq:ode}, starting at $x_i$ at time $t_{n+1}$, computed at time $s=\ddt$, computed with a $\nu$-stage Runge-Kutta (RK) method, and $C^n_i(\underline{a})$ is a proper approximation of the integral of the running cost $f_C$ that appears in \eqref{eq:DP}. Also, 
For a given $\underline{a} \in A^\nu$, $y^n_i(\underline{a})$ is computed via a $\nu$-stage Runge--Kutta (RK) methods as follows:
\begin{equation}\label{eq:RK}
\begin{aligned}
 y^{n}_i(\underline{a}) &= x_i + \ddt\sum_{k=1}^{\nu} b_k K_k(\underline{a}),\\
 K_k(\underline{a}) &= f_D\,(t_n+(1-c_k)\ddt,X_k(\underline{a}),a_k),\\
 X_k(\underline{a}) &= x_i+\ddt\sum_{j=1}^{k-1}A_{kj}K_j(\underline{a}),
\end{aligned} 
\end{equation}
where $b_k$,$c_k$,$A_{kj}$ are the coefficients of the Butcher tableau defining the RK method used to solve \eqref{eq:ode}. For our purposes, we resort to the Forward Euler method for first order schemes, to Heun's method
\begin{equation}\label{eq:heun}
\begin{array}
{c|cc}
0\\
1 & 1\\
\hline
& \nicefrac{1}{2} &\nicefrac{1}{2} 
\end{array}    
\end{equation}
to get second order accuracy, and to RK method with tableau 
\begin{equation}\label{eq:RK3}
\begin{array}
{c|ccc}
0\\
\nicefrac{1}{2} & \nicefrac{1}{2} & \\
1 & -1 & 2 & \\
\hline
& \nicefrac{1}{6} &\nicefrac{2}{3} &\nicefrac{1}{6} 
\end{array}    
\end{equation}
to get third order accuracy.

In order to compute $C^n_i(\underline{a})$, we employ a suitable quadrature rule and replace the function evaluations at the nodes with the numerical solution of \eqref{eq:ode}:
\begin{equation}\label{eq:fcosto}
\begin{aligned}
 &\int_{0}^{\ddt} f_C(t_{n+1}-s,y_{x_i,t_{n+1}}(s),\alpha(s))ds 
 \\
 = &\ddt \sum_k w_k f_C(t_{n+1}- \xi_k \ddt,y_{x_i,t_{n+1}}( \xi_k \ddt),\alpha( \xi_k \ddt)) 
 + \mathcal{O}(\ddt^{Q+1})\\
 = &\ddt \sum_k w_k f_C(t_{n+1}- \xi_k \ddt,\tilde y_{x_i,t_{n+1}}( \xi_k \ddt),\alpha( \xi_k \ddt)) 
 + \mathcal{O}(\ddt^{Q+1}+\ddt^{R+1})
\end{aligned}
\end{equation}
where $\xi_k, w_k$ are the nodes and the weights of the quadrature rule and $\tilde y_{x_i,t_{n+1}}(\xi_k\ddt)$ is a suitable approximation of the numerical solution of $y_{x_i,t_{n+1}}(s)$ for $s\in [0,\ddt]$;
$Q$ is the accuracy of the quadrature rule and $R$ describes the extra error introduced by the approximation of $y_{x_i,t_{n+1}}( \xi_k \ddt)$ with $\tilde{y}_{x_i,t_{n+1}}( \xi_k \ddt)$ within the quadrature rule.
%The quadrature error for the running cost is therefore influenced by both the accuracy of the quadrature rule, and the accuracy of the RK scheme computing $\tilde y$. 
It is tempting to choose a quadrature rule whose nodes coincide with the abscissae $\{c_k\}$ of the RK scheme and replace, in the last approximation of \eqref{eq:fcosto}, $\tilde y_{x_i,t_{n+1}}( \xi_k \ddt)$ with the stage values of the RK scheme $X_{k}$ and $\alpha(\xi_k \ddt)$ with $a_{k}$, for $k=1,\ldots,\nu$.

For the Heun method, the numerical computation of $C^n_i(\underline{a})$ would then be performed with the trapezoidal rule, whose nodes are at $\xi = 0,1$ and for which $\tilde y_{x_i,t_{n+1}}(0)=x_i$ and $\tilde y_{x_i,t_{n+1}}(\ddt)=X_2$ 
are computed with the correct accuracy by the RK scheme \eqref{eq:heun}.

For the third order RK scheme \eqref{eq:RK3}, whose nodes are at $\xi=~0,\nicefrac{1}{2},1$, one would then choose the approximations $\tilde y_{x_i,t_{n+1}}(0)=x_i$, $\tilde y_{x_i,t_{n+1}}( \frac{1}{2}\ddt)= X_2$ and $\tilde y_{x_i,t_{n+1}}(\ddt) = X_3$.
In this special case, even if both $y_{x_i,t_{n+1}}( \xi_2 \ddt)- X_2$ and $y_{x_i,t_{n+1}}( \xi_3 \ddt)-X_3$ are $\Ogrande(\ddx^2)$, these errors cancel each other out when computing the linear combination of the Simpson's quadrature rule, resulting in a third order accuracy for the approximation of $C_i^n(\underline{a})$. For more general RK processes, this issue can be overcome by resorting to the continuous extension \cite{Z86} to compute $\tilde y_{x_i,t_{n+1}}( \xi )$ with the proper accuracy.

% For higher order schemes, the quadrature rule will contain internal nodes, and, in order to compute with the proper accuracy $\tilde y(t_n + \xi_k \ddt)$, we resort to a natural continuous extension (NCE) of the RK method. NCE's were introduced in \cite{Z86}, where also their existence for all RK processes is proven.

% More precisely, given the $\nu$-stage Runge-Kutta process \eqref{eq:RK} of order $p$, its NCE $\tilde y$ of degree $d$ is defined as
% \begin{equation}\label{eq:NCE}
%     \tilde y(t_n+\theta \ddt) \coloneqq x^n_i +\ddt \sum_{s=1}^{\nu} b_s (\theta) K_s, \quad 0\leq \theta\leq 1,
% \end{equation}
% where $b_s(\theta)$, $s=1,...,\nu$, are $\nu$ polynomials of degree $\leq d$, independent of the function $f_D$. For a NCE, the following statements hold:
% \begin{equation}
%     \tilde y(t_n) = x^n_i, \quad \tilde y(t_n+\ddt) = y^n_i
% \end{equation}
% and
% \begin{equation}
%     \underset{t_n \leq t \leq t_n+\ddt}\max | y(t) -\tilde y(t) | = O(\ddt^{d+1}).
% \end{equation}
% For our purposes we will consider the Simpson's rule, the third order SSP Runge-Kutta scheme with tableau...and its NCE with $d=2$, defined as follows:
% \begin{equation}
%  b_1(\theta) = -\frac{1}{2}\theta^2+\frac{2}{3}\theta, \quad
%  b_2(\theta) = \frac{1}{2}\theta^2-\frac{1}{3}\theta, \quad
%  b_3(\theta) = \frac{2}{3}\theta.
% \end{equation}

To compute the solution of \eqref{eq:SL},
we employ tabulation on a coarse grid in $A$, followed by a Nelder-Mead algorithm, adjusted so that no vertex of the simplices can exit the compact set $A$. At each step of the minimization
we need to evaluate the reconstruction at the feet of the characteristics $y^n_i(\underline a)$
for each node $x_i$ and for each discrete control $\underline a$. Thus an efficient and accurate reconstruction operator is crucial to obtain a robust and fast numerical scheme.
Since \eqref{eq:SL} has in general nonsmooth solutions, standard high order interpolation give rise to oscillations in the numerical solution.
\cite{CFR05}, a $\WENO$ reconstruction has been applied to overcome this problem.
The $\CWENO$ approach, instead, is more efficient when many evaluations of the reconstruction in each given cell are required.
For this reason, we propose to apply $\CWENO$ interpolation within our semi-Lagrangian scheme \eqref{eq:SL}.

\section{$\CWENO$ reconstruction}
\label{sec:cweno}
%!TEX root = Cweno_x.tex

We recall here the definition of the $\CWENO$ and $\CWENOZ$ operators, following the presentation of \cite{CPSV:cweno,CSV19:cwenoz}; in particular, despite the present finite-difference setting in place of the finite-volume one, we will still be able to exploit the general theorems for the accuracy of the reconstructions on smooth solutions proven in the aforementioned papers.

For a point $x\in\R^d$, let $\Omega$ be the grid cell containing it and $\stencil_{\Omega}$ be the set of its vertices. In order to achieve better than first order accuracy, we need to consider stencils $\stencil\supset\stencil_\Omega$.
We associate to any $\stencil$ a polynomial $P^{(r)}_{\stencil}(x)\in\mathbb{P}^r_d$ which interpolates the data in $\stencil$, i.e. such that $P^{(r)}_{\stencil}(x_j)=u_j$ $\forall j\in \stencil$.  
Larger symmetric stencils will define highly accurate interpolators, which however would be very oscillatory when the data in $\stencil$ represent a non smooth function. Instead, polynomials associated to smaller stencils, biased in a specific direction, could avoid the discontinuities in the data.

The general idea of $\CWENO$ and $\CWENOZ$ reconstruction is to blend, in a nonlinear and data-dependent fashion, a high order accurate interpolating polynomial with a set of lower order ones to produce an essentially non-oscillatory reconstruction polynomial for the cell $\Omega$.
This can be later evaluated at any point in $x\in\Omega$ with negligible cost.

The nonlinear selection or blending of polynomials relies on oscillation indicators $\OSC[P]$, which are in general scalar quantities associated to a polynomial $P$, designed in such a way that $\OSC[P]\rightarrow 0$ under grid refinement, if $P$ is associated to smooth data, and $\OSC[P] \asymp 1$, in presence of a discontinuity in the stencil $\stencil$. In this work we rely on the classical Jiang--Shu oscillation indicators \cite{JiangShu:96}, suitably modified to accomodate for the regularity of the solution of HJB problems under consideration \cite{JP00,FPT:2020:indHJ}.

We recall now the definition of the  $\CWENO$ and $\CWENOZ$ reconstructions given in \cite{CSV19:cwenoz}.

\begin{Definition}
Given a stencil $\stencil_{\text{opt}}$, including $\stencil_{\Omega}$, let $P_{\text{opt}}\in \mathbb{P}^G_n$ (\textit{optimal polynomial}) be the polynomial of degree $G$, associated to $\stencil_{\text{opt}}$. Further, let $P_1,P_2,...,P_m$ be a set of $m\geq 1$ polynomials of degree $g$ with $g<G$, associated to substencil $\stencil_k$ such that $\stencil_{\Omega}\subset \stencil_k \subset \stencil_{\text{opt}}$ $\forall k=1,...,m$. Let also $\left\{ d_k \right\}_{k=0}^m$ be a set of strictly positive real coefficients such that $\sum_{k=0}^m d_k = 1$.

The $\CWENO$ and $\CWENOZ$ operators compute a reconstruction polynomial 
\begin{equation}
    \begin{aligned}
        P_{\text{rec}}^{CW} &= \CWENO(P_{\text{opt}},P_1,...,P_m) \in \mathbb{P}^G_n,\\
        P_{\text{rec}}^{CWZ} &= \CWENOZ(P_{\text{opt}},P_1,...,P_m) \in \mathbb{P}^G_n,
    \end{aligned}
\end{equation}
as follows:
\begin{enumerate}
    \item First, introduce the polynomial $P_0$ defined as
    \begin{equation}
        P_0(x) = \frac{1}{d_0}\left( P_{\text{opt}}(x) - \sum_{k=1}^m d_k P_k(x) \right) \in \mathbb{P}^G_n;
    \end{equation}
    \item compute suitable regularity indicators
    \begin{equation}
        \OSC_0 = \OSC[P_{\text{opt}}], \quad \OSC_k = \OSC[P_k], \;k\geq 1;
    \end{equation}
    \item compute the nonlinear coefficients $\left\{ \omega_k \right\}_{k=0}^m$ or $\left\{ \omega^Z_k \right\}_{k=0}^m$ as
    \begin{enumerate}
        \item $\CWENO$ operator: for $k=0,...,m$,
        \begin{equation}\label{eq:CWENOomega}
        \alpha_k = \frac{d_k}{(\OSC_k + \epsilon)^l}, \quad \omega_k^{CW} = \frac{\alpha_k}{\sum_{i=0}^m \alpha_i},
        \end{equation}
        \item $\CWENOZ$ operator: for $k=0,...,m$,
        \begin{equation}\label{eq:CWENOZomega}
        \alpha_k^Z = d_k\left(1+\left(\frac{\tau}{\OSC_k + \epsilon}\right)^l\right), \quad \omega_k^{CWZ} = \frac{\alpha_k^Z}{\sum_{i=0}^m \alpha_i^Z},
        \end{equation}
    \end{enumerate}
    where $\epsilon$ is a small positive quantity, $l\geq 1$, and, in the case of $\CWENOZ$, $\tau$ is a global smoothness indicator; 
    \item finally, define the reconstruction polynomial as
    \begin{equation}\label{eq:CWENOrec}
        P_{\text{rec}}^{CW}(x) = \sum_{k=0}^m \omega_k^{CW} P_k(x) \in \mathbb{P}^G_n,
    \end{equation}
    \begin{equation}\label{eq:CWENOZrec}
        P_{\text{rec}}^{CWZ}(x) = \sum_{k=0}^m \omega_k^{CWZ} P_k(x) \in \mathbb{P}^G_n.
    \end{equation}
\end{enumerate}
\end{Definition}

Note that the reconstruction polynomial defined in \eqref{eq:CWENOrec} and in \eqref{eq:CWENOZrec} can be evaluated at any point in the computational cell $\Omega$ at a very low computational cost and this happens because the reconstruction procedure, in both cases, starts with the definition of the linear coefficients $\left\{ d_k \right\}_{k=0}^m$ that do not depend on the reconstruction point. Therefore, the nonlinear coefficients given by \eqref{eq:CWENOomega} and \eqref{eq:CWENOZomega} are computed once per cell and not once per reconstruction point, as in the standard $\WENO$. We also remark that, since every polynomial taking place in the reconstruction procedures is required to satisfy the interpolation constraint on $\stencil_{\Omega}$, then also $P_{\text{rec}}^{CW}$ and $P_{\text{rec}}^{CWZ}$ satisfy the same constraint on the vertices of $\Omega$.

The accuracy and non-oscillatory properties of $\CWENO$ and $\CWENOZ$ schemes are guaranteed by the dependence of their nonlinear weights \eqref{eq:CWENOomega} and \eqref{eq:CWENOZomega} on the regularity indicators $\OSC_k$. On smooth data, the nonlinear weights are driven sufficiently close to the optimal ones, so that $P_{\text{rec}}\approx P_{\text{opt}}$ and the reconstruction reaches the optimal order of accuracy $G+1$. 
On the other hand, when a discontinuity is present in $\stencil_{\text{opt}}$, both $\OSC_0\asymp1$ and at least one $\OSC_{\hat{k}} \asymp 1$ for some $\hat{k}\in \{ 1,...,m \}$.
Then the formulas \eqref{eq:CWENOomega} and \eqref{eq:CWENOZomega} for the nonlinear weights will ensure that $\omega_0 \approx 0$ and 
$\omega_{\hat{k}}\approx 0$ for all $\hat{k}$ such that $P_{\hat{k}}$ would bring oscillations in the reconstruction. The reconstruction polynomial will then be a linear combination of all polynomials of degree $g$ that are not affected by the discontinuity;
the accuracy of the reconstruction thus reduces to $g+1$, but spurious oscillations would be tamed.

\subsection{One spatial dimension}
%In one space dimension, let us consider a spatial domain $[a,b]$ and a uniform grid with points $x_j=a+j\ddx$. 
Let us describe the reconstruction for any point $x$ in the cell $\Omega=[x_j,x_{j+1}]$, so that $\stencil_{\Omega}=\{j,j+1\}$.
We consider $\stencil_{\text{opt}}=\{j-1,j,j+1,j+2\}$ and 
introduce the optimal cubic polynomial $P_{\text{opt}}(x)=Q(x)=\sum_{i=0}^{3}z_i[(x-x_j)/\ddx]^i$ that interpolates the data $u_{j-1},u_j,u_{j+1},u_{j+2}$ at the nodes $x_{j-1},x_{j},x_{j+1},x_{j+2}$. Next, we consider two parabolas $P_L(x)$ and $P_R(x)$ that interpolate only the data $u_{j-1},u_j,u_{j+1}$
and, respectively, $u_j,u_{j+1},u_{j+2}$.

The reconstruction operators thus compute
\begin{equation}
    \begin{aligned}
        P_{\text{rec}}^{CW} &= \CWENO(Q,P_L,P_R) \in \mathbb{P}^3_1,\\
        P_{\text{rec}}^{CWZ} &= \CWENOZ(Q,P_L,P_R) \in \mathbb{P}^3_1.
    \end{aligned}
\end{equation}

For any polynomial, its oscillation indicator is defined as
\begin{equation}
\label{eq:ind:1d}
\OSC[P] = \sum_{\alpha\geq2} \ddx^{2\alpha-3} \int_{x_j}^{x_{j+1}} \left(\frac{d^{(\alpha)}P}{dx^\alpha}\right)^2 dx.
\end{equation}
The above definition is similar to the classical definition of the oscillation indicators for the $\WENO$ reconstruction as given in \cite{JiangShu:96}, except that the first derivative is not included in the sum. This choice is the appropriate one for Hamilton--Jacobi equations, whose solution can be at worst continuous with kinks, see \cite{JP00}, and ensures that at worst $\OSC[P]=\Ogrande(1)$.

Let $P(x)=\sum_{i=0}^{3}z_i[(x-x_j)/\ddx]^i$ be a polynomial of degree up to $3$. Its indicator \eqref{eq:ind:1d} can be written as a quadratic form of its coefficients given by
\begin{equation*}
\OSC[P] = \frac{1}{\ddx^2}(4 z_2^2 + 12z_2 z_3 + 48 z_3^2) 
\end{equation*}
or equivalently, denoting with $\vec{z}$ the vector of coefficients,
\begin{equation}
\OSC[P] = \frac{1}{\ddx^2}\Vec{z}^T M \Vec{z} \text{, with} \quad M = \begin{pmatrix}
0 & 0 & 0 & 0 \\
0 & 0 & 0 & 0 \\
0 & 0 & 4 & 6 \\
0 & 0 & 6 & 48\\
\end{pmatrix}.
\end{equation}
Denoting by $U$ the vector of data $(u_{j-1},u_j,u_{j+1},u_{j+2})^T$, since the coefficients of the polynomial linearly depends on $U$,
we can express the regularity indicators also as a quadratic form on the data being interpolated. 
More precisely, the vector of coefficients $\vec{z}$ can be expressed as $\vec{z} = V^{-1} U$, where $V^{-1}$ is the inverse of the Vandermonde matrix. Thus, in the one-dimensional case, one obtains that the matrices of the quadratic forms expressed in terms of $U$ scale again globally as $1/\ddx^2$:
\begin{equation}\label{eq:ind1dDati}
\OSC[Q] = \frac{1}{\ddx^2} U^T A_Q U, \quad \OSC[P_k] = \frac{1}{\ddx^2} U^T A_k U \quad (k=L,R).
\end{equation}
In our case, we get
\begin{equation*}
A_Q = \begin{pmatrix}
4/3 & -7/2 & 3 & -5/6 \\
-7/2 & 10 & -19/2 & 3 \\
3 & -19/2 & 10 & -7/2 \\
-5/6 & 3 & -7/2 & 4/3 \\
\end{pmatrix},
\end{equation*}
\begin{equation*}
A_L = \begin{pmatrix}
1 & -2 & 1 & 0 \\
-2 & 4 & -2 & 0 \\
1 & -2 & 1 & 0 \\
0 & 0 & 0 & 0 \\
\end{pmatrix}, \hspace{.2cm}
A_R = \begin{pmatrix}
0 & 0 & 0 & 0 \\
0 & 1 & -2 & 1 \\
0 & -2 & 4 & -2 \\
0 & 1 & -2 & 1 \\
\end{pmatrix}.
\end{equation*}

The $\CWENO$ reconstruction applied in the numerical tests of this paper is then defined by choosing linear coefficients $d_L=d_R=\nicefrac{1}{8}$ and $d_0=\nicefrac{3}{4}$, $l=2$ and $\epsilon = \ddx^2$ in the above construction.

For the optimal definition of the $\CWENOZ$ reconstruction, in order to exploit the results of \cite{CSV19:cwenoz}, we need to study the Taylor expansions of the indicators centered at the point $x_0=(x_j+x_{j+1})/2$, namely the center of the reconstruction cell $\Omega$. They are given by
\begin{equation}
    \begin{aligned}
        \OSC[P_L] &= B - u''(x_0)u'''(x_0)\ddx^3 + \frac{1}{4}u'''(x_0)^2\ddx^4 + \Ogrande(\ddx^5),\\
        \OSC[P_R] &= B + u''(x_0)u'''(x_0)\ddx^3  + \frac{1}{4}u'''(x_0)^2\ddx^4 + \Ogrande(\ddx^5),\\
        \OSC[P_L] &= B + \frac{13}{12}u'''(x_0)^2\ddx^4 + \Ogrande(\ddx^5),\\
    \end{aligned}
\end{equation}
where $B=u''(x_0)^2\ddx^2$.
Thus, defining
\begin{equation}\label{eq:tau:1d}
\tau = \bigg| 2\OSC[Q]-\OSC[P_L]-\OSC[P_R] \bigg|,
\end{equation}
all terms up to $\Ogrande(\ddx^3)$ cancel. Note that this is the optimal definition of $\tau$ since it is never possible to cancel all the $\Ogrande(\ddx^4)$ terms in a convex combination of the three indicators. The hypotheses of Theorem $24$ in \cite{CSV19:cwenoz} hold and we can guarantee the optimal order of accuracy of the $\CWENOZ$ reconstruction for $l=2$.

\subsection{Two spatial dimensions}
In two space dimensions, we will not rely on dimensional splitting. Such an approach would in fact nullify the advantages of a fast evaluation of the reconstruction polynomial on a large number of arbitrary points in the reconstruction cell. Our reconstruction will instead generate a polynomial of two variables, like it is done in \cite{SCR:cwenoAMR,SC:19:wb2d} for Cartesian meshes or in general mesh settings \cite{ZhuQiu:2020:triHJ,ADER_CWENO,BGFB:2020}.

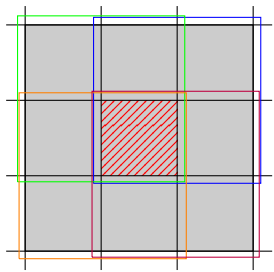
\begin{figure}
\centering
\begin{tikzpicture}
\usetikzlibrary{patterns}
\definecolor{stencil}{gray}{0.8}
\colorlet{cella}{red}

\begin{scope}
\filldraw[fill=stencil](-1.5,-1.5) rectangle (1.5,1.5);
\fill[pattern=north east lines, pattern color=cella] (-.5,-.5) rectangle (.5,.5);
\foreach \x in {-1.5,-0.5,0.5,1.5}
{
  \draw (\x,-1.75) -- (\x,1.75);
  \draw (-1.75,\x) -- (1.75,\x);
}
\draw [blue, xshift=-.6cm,yshift=-.6cm] 
(0,0) 
-- (0,2.2) 
-- (2.2,2.2) 
-- (2.2,0) 
-- (0,0);
\draw [green, xshift=-1.6cm,yshift=-0.58cm] 
(0,0) 
-- (0,2.2) 
-- (2.2,2.2) 
-- (2.2,0) 
-- (0,0);
\draw [orange, xshift=-1.58cm,yshift=-1.6cm] 
(0,0) 
-- (0,2.2) 
-- (2.2,2.2) 
-- (2.2,0) 
-- (0,0);
\draw [purple, xshift=-0.62cm,yshift=-1.58cm] 
(0,0) 
-- (0,2.2) 
-- (2.2,2.2) 
-- (2.2,0) 
-- (0,0);
\end{scope}
\end{tikzpicture} 
\caption{Stencils of the two-dimensional $\CWENO$ and $\CWENOZ$ reconstructions. The red hatched region represents the cell $\Omega_{i,j}$ in which we compute the reconstruction. The vertices of $\stencil_{\text{opt}}$ are enclosed in the grey shaded region, while the stencils for the low degree polynomials are enclosed in the coloured squares: blue, green, orange and purple, respectively for the north-east, north-west, south-west and south-east polynomial.}
\label{fig:CW2dstencils}	
\end{figure}

Let $\Omega$ be the cell of the Cartesian grid containing the reconstruction point $x$.
The setup of the stencils for the reconstruction is illustrated in Fig.~\ref{fig:CW2dstencils}: the cell $\Omega$, with corners $(x_i,y_j)$ and $(x_{i+1},y_{j+1})$, is hatched in red and the reconstruction stencil is shaded in gray.

The reconstruction operators compute 
\begin{equation}
    \begin{aligned}
        P_{\text{rec}}^{CW} &= \CWENO(Q^{(3)}; Q^{(2)}_{\text{ne}}, Q^{(2)}_{\text{se}}, Q^{(2)}_{\text{sw}}, Q^{(2)}_{\text{nw}}) \in \mathbb{P}^3_1 \otimes \mathbb{P}^3_1,\\
        P_{\text{rec}}^{CWZ} &= \CWENOZ(Q^{(3)}; Q^{(2)}_{\text{ne}}, Q^{(2)}_{\text{se}}, Q^{(2)}_{\text{sw}}, Q^{(2)}_{\text{nw}} ) \in \mathbb{P}^3_1 \otimes \mathbb{P}^3_1,
    \end{aligned}
\end{equation}
with the optimal polynomial set to the bicubic polynomial interpolating the 16 data points in the square with corners $(x_{i-1},y_{j-1})$ and $(x_{i+2},y_{j+2})$ enclosed in the grey shaded region of Fig~\ref{fig:CW2dstencils}. The other four polynomials are the biquadratic polynomials that interpolate the four $3\times3$ substencils in the north-east, south-east, south-west and north-west directions, as illustrated in Fig.~\ref{fig:CW2dstencils}, named $\stencil_{\text{ne}}$, $\stencil_{\text{nw}}$, $\stencil_{\text{sw}}$ and $\stencil_{\text{se}}$, respectively. In particular, the blue box with vertices $(x_{i},y_{j})$ and $(x_{i+2},y_{j+2})$ encloses $\stencil_{\text{ne}}$, the green box with vertices $(x_{i-1},y_{j})$ and $(x_{i+1},y_{j+2})$ encloses $\stencil_{\text{nw}}$, the orange box with vertices $(x_{i-1},y_{j-1})$ and $(x_{i+1},y_{j+1})$ encloses $\stencil_{\text{sw}}$ and the purple box with vertices $(x_{i},y_{j-1})$ and $(x_{i+2},y_{j+1})$ encloses $\stencil_{\text{se}}$. The basis for bicubic and biquadratic polynomials are defined by tensorization and are thus respectively given by
\begin{align*}
    \mathcal{B}_3 &= \{ 1,\hat{x},\hat{y},\hat{x}^2,\hat{x}\hat{y},\hat{y}^2,\hat{x}^3,\hat{x}^2\hat{y},\hat{x}\hat{y}^2,\hat{y}^3,\hat{x}^3\hat{y},\hat{x}^2\hat{y}^2,\hat{x}\hat{y}^3,\hat{x}^3\hat{y}^2,\hat{x}^2\hat{y}^3,\hat{x}^3\hat{y}^3\},\\
    \mathcal{B}_2 &= \{ 1,\hat{x},\hat{y},\hat{x}^2,\hat{x}\hat{y},\hat{y}^2\},
\end{align*}
where $\hat{x}=(x-x_i)/\ddx$ and $\hat{y}=(y-y_j)/\ddy$.

The oscillation indicators are defined as
\begin{equation}\label{eq:IND:2d}
\OSC[P] 
= 
\sum_{|\vec{\alpha}|\geq2} \Delta^{2|\vec{\alpha}|-4} 
\int_{\Omega_{i,j}}
\left[
  \frac{\partial P}{\partial x^{\alpha_1}\partial y^{\alpha_2}}
\right]^2
\dx\dy
\end{equation}
where $\Delta$ is the diameter of the cell $\Omega_{i,j}$.
These indicators are inspired to the classical ones of \cite{HuShu:99}, but the powers of $\Delta$ appearing in \eqref{eq:IND:2d} have been adjusted for the expected regularity of the solution to Hamilton--Jacobi equations.
In particular the scaling $\Delta^{2|\vec{\alpha}|-4}$ 
is a generalization to two space dimensions of the choice of \cite{JP00};
it is introduced so that $\OSC[P]\asymp1$ for functions with discontinous first derivative and $\OSC[P]=\Ogrande(\Delta^2)$ on regular solutions. A detailed discussion of this scaling can be found in \cite{FPT:2020:indHJ}.

Denoting with $\vec{z}$ the vector of coefficients along the basis $\mathcal{B}_3$ of a generic polynomial in two variables of degree up to 3, and considering a uniform grid of size $\ddx$, its oscillation indicator can be expressed as
\begin{equation}\label{eq:ind2dCoeff}
\OSC[P] = \frac{1}{\ddx^2}\Vec{z}^T M \Vec{z}
\end{equation}
for a $16\times 16$ matrix $M$%, which is reported in the Appendix.
We point out that in the case of Cartesian but uneven grids,
there would still be a global factor inversely proportional to the local grid size, but matrix $M$ would also depend on the aspect ratios and on the ratios of sizes of nearby cells.

Consider now the stencil of the reconstruction given by $\{(x_{i-1+k},y_{j-1+l})\}_{\tiny{\substack{l=0,\ldots,3\\k=0,\ldots,3}}}$ and denote by $U$ the corresponding vector of values $\{(u_{i-1+k},y_{j-1+l})\}_{\tiny{\substack{l=0,\ldots,3\\k=0,\ldots,3}}}$ with a prescribed ordering, e.g. lexicographic. Thanks to the linear relation among the coefficients $\vec{z}$ and the data vector $U$, the regularity indicators can be expressed in the form
\begin{equation}\label{eq:ind2dDati}
\OSC[Q] = \frac{1}{\ddx^2}U^T A_Q U, \quad \OSC[Q_k] = \frac{1}{\ddx^2}U^T A_k U \quad (k=\text{ne},\text{se},\text{sw},\text{nw}).
\end{equation}
The matrices $A_Q$ and $A_k$ will in general depend on the local aspect ratio of cells and on the neighbours size ratio. 
%Their entries are reported in the Appendix for the case of a uniform Cartesian grid.

The $\CWENO$ reconstruction is then defined by choosing linear coefficients $d_k=1/16$ and $d_0=1-\sum_k d_k$ with $k \in \{\text{ne},\text{se},\text{sw},\text{nw}\}$, $l=2$ and $\epsilon = \ddx^2$.

The $\CWENOZ$ reconstruction is computed choosing $l=2$ and
\begin{equation}\label{eq:tau:2d}
\tau = \left| 
4\OSC[Q^{(3)}]
-\OSC[Q^{(2)}_{\text{ne}}]
-\OSC[Q^{(2)}_{\text{se}}]
-\OSC[Q^{(2)}_{\text{sw}}]
-\OSC[Q^{(2)}_{\text{nw}}]
\right|
\end{equation}
in order to guarantee the optimal order of accuracy of the $\CWENOZ$ reconstruction, see \cite[Theorem 24]{CSV19:cwenoz}. In fact, computing the Taylor expansions of the indicators centered at the center of the cell $\Omega$, one gets, for smooth enough data,
\begin{equation}
    \begin{aligned}
        \OSC[Q^{(2)}_{\text{ne}}] &= B + 2u_{xx}u_{xxx}\ddx^3 + 2u_{yy}u_{yyy}\ddx^3 + \Ogrande(\ddx^4),\\
        \OSC[Q^{(2)}_{\text{nw}}] &= B - 2u_{xx}u_{xxx}\ddx^3 + 2u_{yy}u_{yyy}\ddx^3 + \Ogrande(\ddx^4),\\
        \OSC[Q^{(2)}_{\text{sw}}] &= B - 2u_{xx}u_{xxx}\ddx^3 - 2u_{yy}u_{yyy}\ddx^3 + \Ogrande(\ddx^4),\\
        \OSC[Q^{(2)}_{\text{se}}] &= B + 2u_{xx}u_{xxx}\ddx^3 - 2u_{yy}u_{yyy}\ddx^3 + \Ogrande(\ddx^4),\\
        \OSC[Q^{(3)}] &= B + \Ogrande(\ddx^4)\\
    \end{aligned}
\end{equation}
where $B=(u_{xx}^2+u_{xy}^2+u_{yy}^2)\ddx^2$. Therefore, the coefficients chosen for the definition of $\tau$ in \eqref{eq:tau:2d}, cancel all terms up to $\Ogrande(\ddx^3)$, thus ensuring that the hypotheses of \cite[Theorem 24]{CSV19:cwenoz} hold.

\section{Convergence} 
\label{sec:convergence}
%!TEX root = Cweno.tex

%\todo[inline]{Richiamare qui la definizione dell' ``operatore \CWENO'' da MathComp e \CWENOZ\ da SINUM e riscrivere entrambe le ricostruzioni con questo formalismo}

The convergence analysis will be carried out in one space dimension, and follow the guidelines of \cite{F01,CFR05}, which will be briefly reviewed here. Assume that equation \refp{eq:HJ} is posed on $\R$ in the simplified form
\begin{equation}\label{eq:HJ_simpl}
\begin{cases}
v_t(t,x)+H(D v(t,x))=0 & t,x \in(0,T)\times \R,\\
v(0,x)=v_0(x), &x\in \R. 
\end{cases}
\end{equation}
%and discretized on an infinite uniform grid with nodes $x_j=j\ddx$ $(j=0,\pm1,\pm2,\ldots)$. 
Also assume that the Hamiltonian
function $H$ is a $W^{2,\infty}$ function and that
\begin{equation}\label{eq:hamilt}
H''(p) \ge m_H > 0.
\end{equation}
The convexity assumption \eqref{eq:hamilt} implies, by the Fenchel duality formula, that $H(\cdot)$ can always be written in the form
\[
H(Dv(t,x)) = \sup_{a \in \mathbb{R}} \{a\, Dv(t,x) - H^*(a)\},
\]
where \(H^*\) denotes the Legendre transform of \(H\). Moreover, the restriction of the supremum can be limited to a suitable compact set \(A\), so that problem \eqref{eq:HJ_simpl} can be recast in the case of a Hamiltonian of the form \eqref{eq:Hamiltonian}, by setting \(f_D(t,x,a) = -a\) and \(f_C(t,x,a) = H^*(a)\).

The key assumption is that for any Lipschitz continuous function $v(x)$, once defined the sequence $V=\{v_j\}_{j}=\{v(x_j)\}_{j}$, the interpolation operator $R[V]$ satisfies $R[V](x_k)=v(x_k)$ and, for some constant $C<1$:
\begin{equation}\label{eq:hyp}
|R[V](x)-R_1[V](x)|\le C \max_{x_k\in U(x)}
|v_{k+1}-2v_k+v_{k-1}|
\end{equation}
where by $R_1[V]$ we denote the $\mathbb{P}_1$ (i.e., piecewise linear) interpolation on the sequence $V$, and by $U(x)=(x-h_-\ddx,x+h_+\ddx)$ the interval containing the stencil of the reconstruction $R[V](x)$. For example, in our case the reconstruction is performed taking two nodes on the left and two nodes on the right of the point $x$, so that $h_-=h_+=2$.

Under such assumptions it is possible (see \cite{F01}) to prove the following
\begin{Theorem}\label{thm:conv_thm}
Consider the scheme \eqref{eq:SL} applied to equation \eqref{eq:HJ_simpl}, and assume that \eqref{eq:hamilt} and \eqref{eq:hyp} hold, that $\ddx=\Ogrande(\ddt^2)$ and that
$v_0$ is Lipschitz continuous. Then, the numerical solution $U^n=\{u_j^n\}_{j}$ (with $u_j^n$
defined by \eqref{eq:SL}) satisfies
$$
\|R[U^n]-v(n\ddt)\|_\infty\to 0
$$
(where $v$ is the solution of \eqref{eq:HJ_simpl}) for $0\le n\le T/\ddt$, as $\ddt\to 0$.
\end{Theorem}

It is proved in \cite{F01} that condition \refp{eq:hyp} is satisfied for Lagrange reconstructions up to the fifth order, if the reconstruction stencil includes the interval $[x_j,x_{j+1}]$. Since $\WENO$
interpolation is performed by taking a convex combination of polynomials which satisfy \refp{eq:hyp} up the degree 5 for the partial polynomials (and therefore, up to the degree 9 for the global interpolant), this fact is used in \cite{CFR05} to obtain convergence via Theorem \ref{thm:conv_thm}. This essentially boils down to proving that the all the function appearing in the convex combination defining $R$ satisfy \eqref{eq:hyp}, and we plan to apply the same principle here.

In the case of $\CWENO$, once set $d_L=d_R=d$ and recast $P_0$ as
\[
P_0(x) = \frac{1}{1-2d}(Q(x)-dP_L(x)-dP_R(x)),
\]
we estimate the left-hand side of \eqref{eq:hyp} as
\begin{eqnarray*}
|R-R_1| & = & |\omega_0P_0+\omega_LP_L+\omega_RP_R-R_1| \\
& = & |\omega_0(P_0-R_1)+\omega_L(P_L-R_1)+\omega_R(P_R-R_1)| \\
& \le & \max \{|P_0-R_1|,|P_L-R_1|,|P_R-R_1| \},
\end{eqnarray*}
where the identity $\omega_0+\omega_L+\omega_R=1$ has been used twice.
As it has been proved in \cite{F01}, both $P_L$ and $P_R$ satisfy \eqref{eq:hyp}. It remains then to check that the same property holds for $P_0$. We have:
\begin{eqnarray*}
|P_0-R_1| & = & \left|\frac{1}{1-2d}(Q-dP_L-dP_R)-R_1\right| \\
& = & \left|\frac{1}{1-2d}(Q-dP_L-dP_R+(2d-1)R_1)\right| \\
& = & \left|\frac{1}{1-2d}\big((Q-R_1)-d(P_L-R_1)-d(P_R-R_1)\big)\right| \\
& \le & \frac{1}{1-2d}|Q-R_1|+\frac{d}{1-2d}|P_L-R_1|+\frac{d}{1-2d}|P_R-R_1|.
\end{eqnarray*}
Denoting now by $C_r$ the constant appearing in \eqref{eq:hyp} for the specific case of an interpolation of degree $r$, we finally obtain an estimate in the form
\begin{equation*}
|P_0-R_1| \le \left(\frac{1}{1-2d}C_3+\frac{2d}{1-2d}C_2\right) \max_{x_k\in U(x)}
|v_{k+1}-2v_k+v_{k-1}|.
\end{equation*}
As proved in \cite{F01}, $C_2=1/8$ and $C_3\approx 0.2533$; then, it turns out that \eqref{eq:hyp} is satisfied for $d\lesssim 0.332$, and therefore, with this additional condition, all the assumptions of Theorem \ref{thm:conv_thm} are satisfied, and the scheme \eqref{eq:SL} converges to the viscosity solution of \eqref{eq:HJ_simpl}.

We point out that our choice of $d_0=3/4$, i.e. $d=d_L=d_R=1/8$ in the $\CWENO$ and $\CWENOZ$ reconstructions employed in this paper fully satisfies the above requirement.

\section{Numerical tests}
\label{sec:tests}

Here, we present some numerical tests in order to assess the performance of the schemes proposed in this paper. In particular, we focus on the expected accuracy of the schemes, considering the $L^1$-norm of the error, and the computational times, reported in seconds.
In particular we will compare the SL schemes using the $\CWENO$ and $\CWENOZ$ reconstructions of \S\ref{sec:cweno} and those employing the dimensionally-split $\WENO$ approach of \cite{CFR05}.

In order to account for bounded domains, we consider either periodic boundary conditions or extrapolation technique.
 It is noteworthy that in all the numerical tests where extrapolation is applied, the solution has a compact support, or the boundary of the numerical domain is entirely of outflow type.  In both cases, these conditions ensure no loss in accuracy in the numerical treatment of the boundary. In all tests, unless otherwise stated, we employ an extrapolation technique.
% Periodic boundary conditions are imposed by treating the domain as a torus and periodically mapping the characteristics $y_i^n(\underline{a})$ within the torus.
For more general boundary conditions, we refer to \cite{CCDS} for approximations of Neumann-type boundary conditions and to \cite{BCCF} for a second order accurate treatment of Dirichlet boundary conditions.

Although the convergence analysis requires a relation $\ddx = \Ogrande(\ddt^2)$ between the time and space steps, in practice, smaller time steps are also allowed. Heuristically, assuming exact minimization, and supposing a third order time discretization is used, then the consistency error is given by
\[
\T(\ddt,\ddx) = \Ogrande\left(\ddt^3 + \frac{\ddx^q}{\ddt}\right),
\]
where $q$ is the order of the interpolation in space. For smooth solutions, the choice $\ddt = \Ogrande\left(\ddx^{q/4}\right)$ optimizes this error and leads to a numerical convergence rate of $3q/4$.

In particular, choosing $q = 4$ (i.e., a cubic reconstruction), we obtain
\[
\T(\ddt,\ddx) = \Ogrande\left(\ddt^3 + \frac{\ddx^4}{\ddt}\right),
\]
which results in an overall order of $3$ for the choice $\ddt = \Ogrande(\ddx)$. In the particular case where the Hamiltonian is independent of $(t,x)$, the characteristics solving problem \eqref{eq:ode} are affine, and even a first order time discretization computes them exactly, reducing the consistency error to
\[
\T(\ddt,\ddx) = \Ogrande\left( \frac{\ddx^q}{\ddt}\right).
\]
In this case, the larger the time step, the smaller the consistency errors.

 In what follows $N$ indicates the number of grid points per direction, so that $\ddx\sim1/N$. The tests have been performed on the cluster Galileo 100 hosted at CINECA\footnote{https://www.hpc.cineca.it/systems/hardware/galileo100/}.
\paragraph{Test 1: Passive advection.}
\begin{figure}
\centering
\epsfig{file=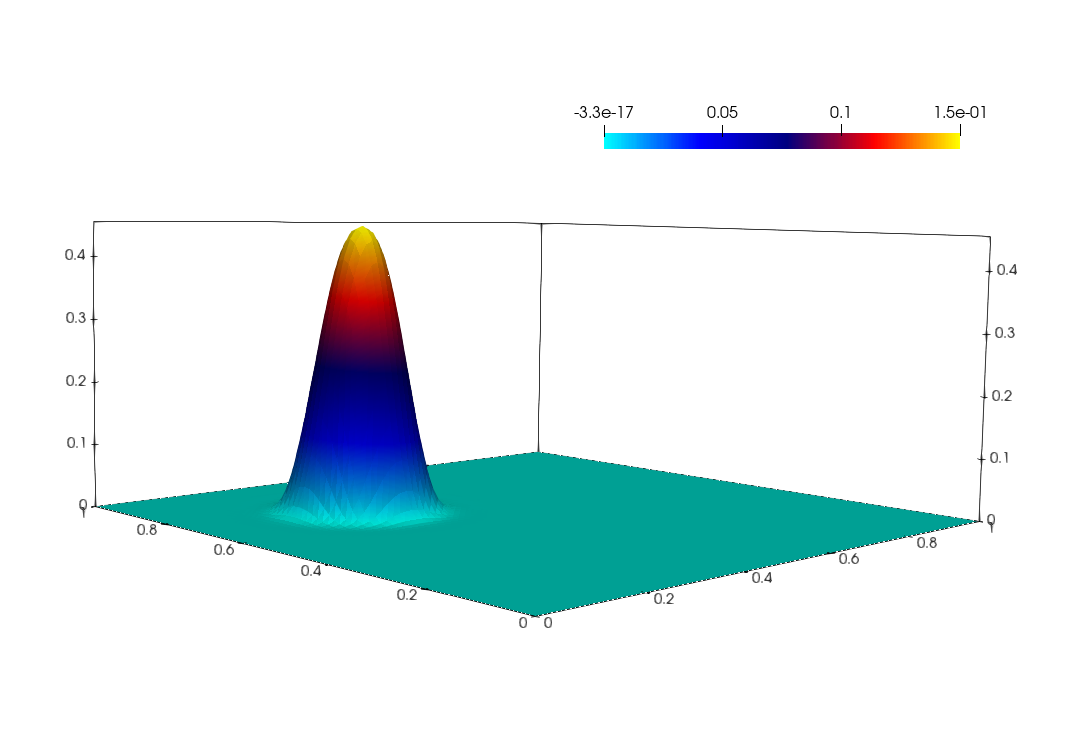,width=0.45\linewidth}
\epsfig{file=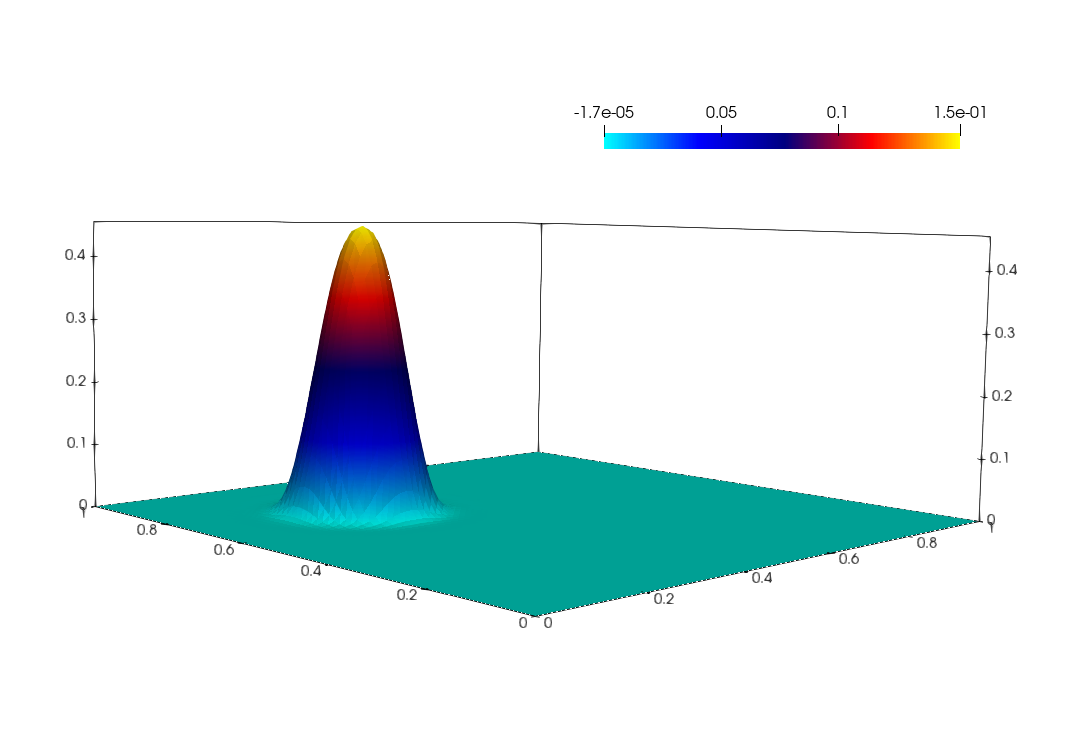,width=0.45\linewidth}
\caption{Initial condition (left) and numerical solution at $T=1$ (right) for Test 1 computed with $\CWENOZ$ reconstruction and $N=81$.}\label{fig:test1}
\end{figure}

\begin{figure}
    \centering
    \includegraphics[width=0.5\linewidth]{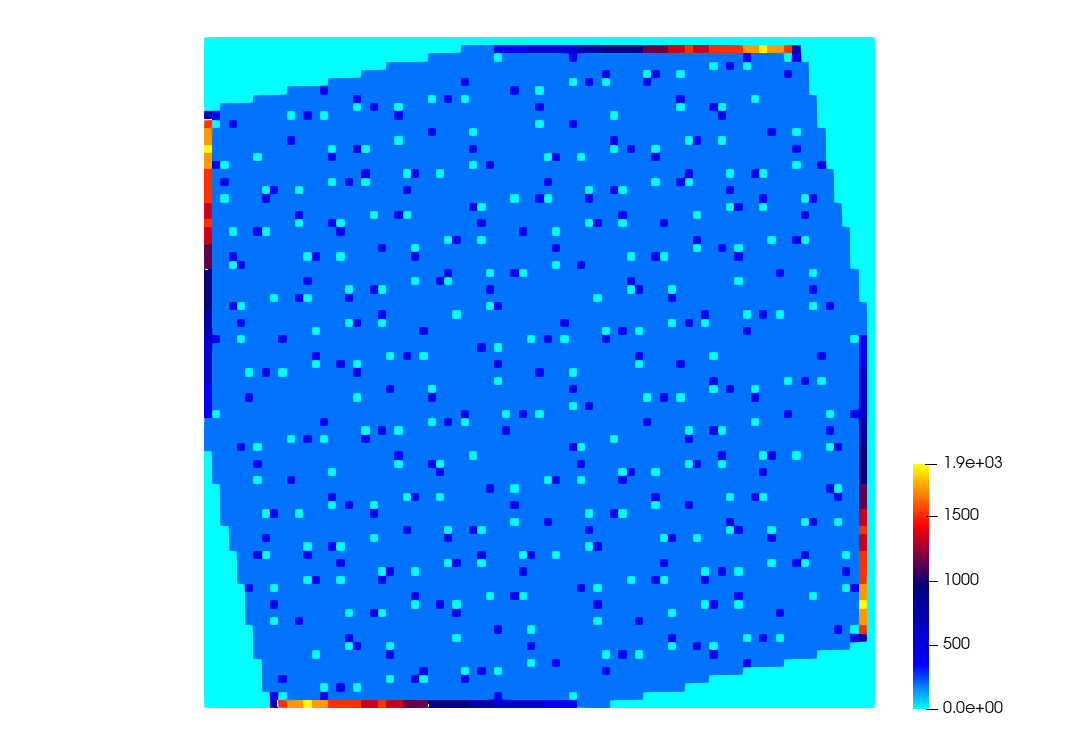}
    \caption{Reconstructions count for Test 1 in the final time step on a grid of $81\times81$ points.}
    \label{fig:test1recCount}
% \todo[inline]{Mi verrebbe voglia di provare a togliere le celle al bordo... In paraview c'è un'icona ``extract subset'': terza riga di icone sotto i menu, stessa fila della calcolatrice e delle slice. Se togliamo quelle celle, la scala colori viene più significativa?}
\end{figure}

\begin{table}
\begin{center}
\begin{tabular}{|c|c|c|c|c|c|c|}
\hline &\multicolumn{2}{ |c|}{$\WENO$} &\multicolumn{2}{|c|}{$\CWENO$} &\multicolumn{2}{|c|}{$\CWENOZ$} \\
\hline $N$ & $L^{1}$-err & $L^1$-ord & $L^{1}$-err & $L^1$-ord & $L^{1}$-err & $L^1$-ord\\ \hline \hline
$21\times21$ & $2.97e-03$ & & $2.98e-03$ & & $2.98e-03$ & \\ \hline 
$41\times41$ & $5.46e-04$ & $2.44$ & $4.88e-04$ & $2.61$ & $4.84e-04$ & $2.62$\\ \hline 
$81\times81$ & $1.24e-04$ & $2.13$ & $8.18e-05$ & $2.57$ & $7.85e-05$ & $2.62$\\ \hline 
$161\times161$ & $2.34e-05$ & $2.41$ & $1.22e-05$ & $2.75$ & $1.12e-05$ & $2.80$\\ \hline 
$321\times321$ & $3.92e-06$ & $2.58$ & $1.72e-06$ & $2.82$ & $1.49e-06$ & $2.92$\\ \hline 
$641\times641$ & $6.11e-07$ & $2.68$ & $2.27e-07$ & $2.92$ & $1.79e-07$ & $3.05$\\ \hline 
$1281\times1281$ & $8.87e-08$ & $2.78$ & $2.92e-08$ & $2.96$ & $2.06e-08$ & $3.12$\\ \hline 
  \end{tabular}
\end{center}
\caption{ Errors at time $T=1$ for Test 1, $\WENO$, $\CWENO$ and $\CWENOZ$  schemes.\label{tab:test1err} }
\end{table}

% CINECA
\begin{table}
\begin{center}
\begin{tabular}{|c|c|c|c|c|c|}
\hline & $\WENO$ & \multicolumn{2}{ |c|}{$\CWENO$} & \multicolumn{2}{ |c|}{$\CWENOZ$}\\
\hline grid & CPU time & CPU time & $\%$ gain & CPU time & $\%$ gain \\ \hline
\hline
$21\times21$ & $1.05e+00$ & $7.22e-01$ & $31.21$ & $1.01e-01$ & $90.36$\\ \hline 
$41\times41$ & $1.10e+00$ & $1.06e+00$ & $3.39$ & $7.20e-01$ & $34.25$\\ \hline 
$81\times81$ & $7.84e+00$ & $5.42e+00$ & $30.87$ & $5.46e+00$ & $30.39$\\ \hline 
$161\times161$ & $6.19e+01$ & $4.36e+01$ & $29.51$ & $4.37e+01$ & $29.35$\\ \hline 
$321\times321$ & $4.86e+02$ & $3.38e+02$ & $30.47$ & $3.36e+02$ & $30.95$\\ \hline 
$641\times641$ & $3.87e+03$ & $2.70e+03$ & $30.31$ & $2.67e+03$ & $30.99$\\ \hline 
$1281\times1281$ & $3.09e+04$ & $2.12e+04$ & $31.25$ & $2.12e+04$ & $31.32$\\ \hline 
  \end{tabular}
\end{center}
\caption{CPU times for Test 1, $\WENO$ and $\CWENO$ schemes.\label{tab:test1cpu} }
\end{table}
As a first benchmark problem, we consider the uniform rotation of a scalar field in two space dimensions, expressed by the linear transport equation
\[
v_t - f_D \cdot Dv = 0,
\]
which is included in our setting by choosing in \eqref{eq:Hamiltonian} $A = \emptyset$ and $f_C = 0$.
The speed of propagation is given by the vector field \(f_D = (-2\pi (x_2-0.5) , 2\pi (x_1-0.5))\) in the square domain \([0,1]\times [0,1]\).
We choose as initial condition the globally $\mathcal{C}^2$ function
$$v(0,x)=v_0(r)= M \left(1 + \frac{r^3}{R^3}\left(-1 + 3\,\frac{r-R}{R}\left(1 - 2\,\frac{r-R}{R}\right)\right) \right)$$
with $r = |x - (0.3,0.7)|$ and parameters $M=0.15$, $R=0.15$.
Numerical solutions are computed at the final time $T=1$, using a relationship of $\ddt= 3\ddx$ among the discretization steps. Since the advection is rigid, the exact solution at final time coincides with the initial condition.
We point out that the characteristic lines are not affine, therefore we use the third order RK scheme and a third order reconstruction scheme as described in section \S\ref{sec:numericalScheme} and \S\ref{sec:cweno}, respectively. Plots of the initial data and of the numerical solution are shown in Fig.~\ref{fig:test1}.

As expected, since the transport is rigid and the initial data is globally $\mathcal{C}^2$, all the schemes ($\WENO$, $\CWENO$ and $\CWENOZ$) yield solutions converging with the same order, see Table \ref{tab:test1err}. However, central reconstructions are more accurate with respect to the traditional $\WENO$ one; in particular, $\CWENOZ$ achieves the best results: its errors are about one-half than those of $\WENO$. The fundamental difference among the schemes is in their computational costs, which are shown in Table \ref{tab:test1cpu}. The high number of reconstructions performed in each cell (see Fig.~\ref{fig:test1recCount}) allows the $\CWENO$ and $\CWENOZ$ to save about the $30\%$ of time. This happens because Central $\WENO$ reconstructions compute their nonlinear coefficients only once per cell instead of once per reconstruction point. Looking at Fig.~\ref{fig:test1recCount}, note in particular that the high number of reconstructions performed near the boundary is related to the boundedness of the domain: in fact, except for the case of periodic boundary conditions, each time a foot of a characteristic falls outside the domain, it is projected back in the first cell near the boundary in order to compute the reconstruction.  
%%%%%%%%%%%%%%%%%%%%%%%%%%%%%%%%%%%%
\paragraph{Test 2: One--dimensional semi-concave data.}

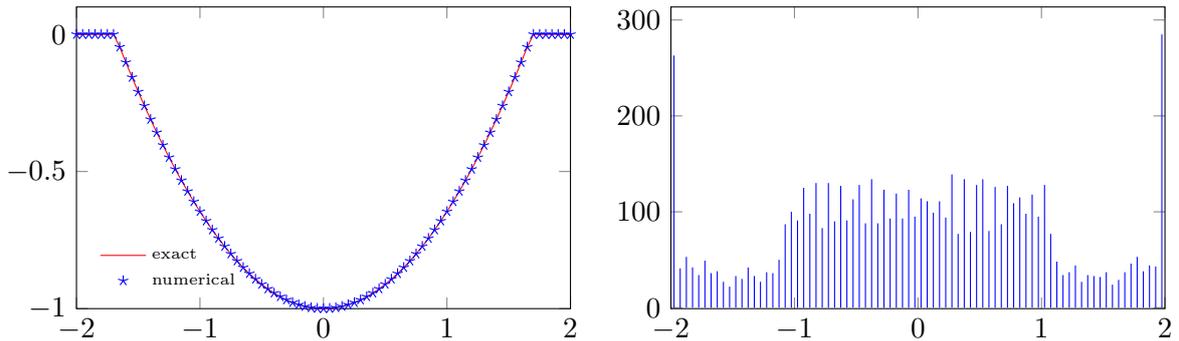
\begin{figure}
\centering
%\documentclass{standalone}
%\usepackage{pgfplots}
%\begin{document}
\begin{tikzpicture}
\pgfplotstableread{%
step variation
-2	0
-1.95000000000000	0
-1.90000000000000	0
-1.85000000000000	0
-1.80000000000000	0
-1.75000000000000	0
-1.70000000000000	0
-1.65000000000000	-0.0461386364631470
-1.60000000000000	-0.102127534392356
-1.55000000000000	-0.156523020987177
-1.50000000000000	-0.209310739472492
-1.45000000000000	-0.260476124675947
-1.40000000000000	-0.310004887931116
-1.35000000000000	-0.357885721692143
-1.30000000000000	-0.404107015502213
-1.25000000000000	-0.448657323132160
-1.20000000000000	-0.491527788613987
-1.15000000000000	-0.532709048346862
-1.10000000000000	-0.572191683403540
-1.05000000000000	-0.609968937768554
-1	-0.646033239406251
-0.950000000000000	-0.680376582593900
-0.900000000000000	-0.712994081442151
-0.850000000000000	-0.743879705575160
-0.800000000000000	-0.773026388590530
-0.750000000000000	-0.800430720995839
-0.700000000000000	-0.826088147176168
-0.650000000000000	-0.849992934766656
-0.600000000000000	-0.872142245770550
-0.550000000000000	-0.892532754759479
-0.500000000000000	-0.911160291296550
-0.450000000000000	-0.928022375812199
-0.400000000000000	-0.943116591780622
-0.350000000000000	-0.956440386849748
-0.300000000000000	-0.967991680078574
-0.250000000000000	-0.977768597245106
-0.200000000000000	-0.985770213446637
-0.150000000000000	-0.991995043657628
-0.0999999999999999	-0.996441395027622
-0.0499999999999998	-0.999109912189189
0	-1
0.0500000000000003	-0.999109912189189
0.100000000000000	-0.996441395027622
0.150000000000000	-0.991995043657628
0.200000000000000	-0.985770213446637
0.250000000000000	-0.977768597245106
0.300000000000000	-0.967991680078574
0.350000000000000	-0.956440386849748
0.400000000000000	-0.943116591780622
0.450000000000000	-0.928022375812199
0.500000000000000	-0.911160291296549
0.550000000000000	-0.892532754759479
0.600000000000000	-0.872142245770550
0.650000000000000	-0.849992934766655
0.700000000000000	-0.826088147176168
0.750000000000000	-0.800430720995839
0.800000000000000	-0.773026388590530
0.850000000000000	-0.743879705575160
0.900000000000000	-0.712994081442151
0.950000000000000	-0.680376582593900
1	-0.646033239406251
1.05000000000000	-0.609968937768553
1.10000000000000	-0.572191683403541
1.15000000000000	-0.532709048346862
1.20000000000000	-0.491527788613987
1.25000000000000	-0.448657323132160
1.30000000000000	-0.404107015502212
1.35000000000000	-0.357885721692142
1.40000000000000	-0.310004887931116
1.45000000000000	-0.260476124675947
1.50000000000000	-0.209310739472492
1.55000000000000	-0.156523020987176
1.60000000000000	-0.102127534392356
1.65000000000000	-0.0461386364631468
1.70000000000000	0
1.75000000000000	0
1.80000000000000	0
1.85000000000000	0
1.90000000000000	0
1.95000000000000	0
2	0
}{\solTestDue}

\pgfplotstableread{%
step variation
-2	0
-1.95000000000000	0
-1.90000000000000	0
-1.85000000000000	0
-1.80000000000000	0
-1.75000000000000	0
-1.70000000000000	0
-1.65000000000000	-0.0461391983921223
-1.60000000000000	-0.102127689845286
-1.55000000000000	-0.156523612136762
-1.50000000000000	-0.209311327727356
-1.45000000000000	-0.260476409066846
-1.40000000000000	-0.310005520069327
-1.35000000000000	-0.357886313866546
-1.30000000000000	-0.404107344235208
-1.25000000000000	-0.448657988521747
-1.20000000000000	-0.491528380287755
-1.15000000000000	-0.532709350168362
-1.10000000000000	-0.572192373869347
-1.05000000000000	-0.609969526133170
-1	-0.646033439995017
-0.950000000000000	-0.680377270579003
-0.900000000000000	-0.712994662869692
-0.850000000000000	-0.743879722970913
-0.800000000000000	-0.773026992442273
-0.750000000000000	-0.800431425361853
-0.700000000000000	-0.826088367815658
-0.650000000000000	-0.849993539560519
-0.600000000000000	-0.872143017639355
-0.550000000000000	-0.892533221761983
-0.500000000000000	-0.911160901288035
-0.450000000000000	-0.928023123674816
-0.400000000000000	-0.943117264268651
-0.350000000000000	-0.956440997338332
-0.300000000000000	-0.967992288261185
-0.250000000000000	-0.977769386789223
-0.200000000000000	-0.985770821330590
-0.150000000000000	-0.991995394195023
-0.0999999999999999	-0.996442177758893
-0.0499999999999998	-0.999110511518164
0	-1
0.0500000000000003	-0.999110511518059
0.100000000000000	-0.996442177758920
0.150000000000000	-0.991995394195039
0.200000000000000	-0.985770821330628
0.250000000000000	-0.977769386789173
0.300000000000000	-0.967992288261118
0.350000000000000	-0.956440997338376
0.400000000000000	-0.943117264268725
0.450000000000000	-0.928023123674800
0.500000000000000	-0.911160901288021
0.550000000000000	-0.892533221761975
0.600000000000000	-0.872143017639364
0.650000000000000	-0.849993539560427
0.700000000000000	-0.826088367815675
0.750000000000000	-0.800431425361890
0.800000000000000	-0.773026992442345
0.850000000000000	-0.743879722970890
0.900000000000000	-0.712994662869679
0.950000000000000	-0.680377270579018
1	-0.646033439994928
1.05000000000000	-0.609969526133222
1.10000000000000	-0.572192373869481
1.15000000000000	-0.532709350168701
1.20000000000000	-0.491528380288649
1.25000000000000	-0.448657988535556
1.30000000000000	-0.404107344258909
1.35000000000000	-0.357886313908783
1.40000000000000	-0.310005520148027
1.45000000000000	-0.260476409220565
1.50000000000000	-0.209311328043118
1.55000000000000	-0.156523612822383
1.60000000000000	-0.102127691427720
1.65000000000000	-0.0461392022953133
1.70000000000000	0
1.75000000000000	0
1.80000000000000	0
1.85000000000000	0
1.90000000000000	0
1.95000000000000	0
2	0
}{\solTestDueExa}

\begin{axis}[%
  width=6.5cm,height=4cm,scale only axis,
  xmin=-2,xmax=2,ymin=-1,
  legend columns=1, 
  legend cell align={left},
  legend style={
        fill=none,
  	legend pos=south west,
  	draw=none,
        font=\tiny,
  	%cells={anchor=east}
  }
]
\addplot[color=red,mark=none] table[x=step,y=variation] {\solTestDueExa};
\addlegendentry{exact};
\addplot[blue,mark=star,only marks] table[x=step,y=variation] {\solTestDue};
\addlegendentry{numerical};

%\draw[dashed,help lines] (axis cs:3.5,0.01) -- (axis cs:3.5,.12) ;
%\draw[dashed,help lines] (axis cs:8.5,0.01) -- (axis cs:8.5,.12) ;
%\draw[stealth-stealth,help lines] (axis cs:0.5,0.02) -- node[pos=0.5,below,black]{\tiny $m=50$} (axis cs:3.5,.02) ;
%\draw[stealth-stealth,help lines] (axis cs:3.5,0.02) -- node[pos=0.5,below,black]{\tiny $m=71$} (axis cs:8.5,.02) ;
%\draw[stealth-stealth,help lines] (axis cs:8.5,0.02) -- node[pos=0.5,below,black]{\tiny $m=100$} (axis cs:20.5,.02) ;

\end{axis}
\end{tikzpicture}
%\end{document} 
%\documentclass{standalone}
%\usepackage{pgfplots}
%\begin{document}
\begin{tikzpicture}
\pgfplotstableread{%
step variation
-1.97500000000000	263
-1.92500000000000	41
-1.87500000000000	53
-1.82500000000000	42
-1.77500000000000	34
-1.72500000000000	49
-1.67500000000000	36
-1.62500000000000	38
-1.57500000000000	27
-1.52500000000000	22
-1.47500000000000	33
-1.42500000000000	30
-1.37500000000000	42
-1.32500000000000	33
-1.27500000000000	27
-1.22500000000000	37
-1.17500000000000	36
-1.12500000000000	50
-1.07500000000000	87
-1.02500000000000	100
-0.975000000000000	91
-0.925000000000000	125
-0.875000000000000	98
-0.825000000000000	130
-0.775000000000000	83
-0.725000000000000	130
-0.675000000000000	90
-0.625000000000000	127
-0.575000000000000	91
-0.525000000000000	113
-0.475000000000000	128
-0.425000000000000	88
-0.375000000000000	134
-0.325000000000000	88
-0.275000000000000	123
-0.225000000000000	93
-0.175000000000000	119
-0.125000000000000	93
-0.0749999999999998	123
-0.0249999999999999	95
0.0250000000000001	114
0.0750000000000002	111
0.125000000000000	99
0.175000000000000	111
0.225000000000000	94
0.275000000000000	139
0.325000000000000	77
0.375000000000000	134
0.425000000000000	79
0.475000000000000	128
0.525000000000000	134
0.575000000000000	80
0.625000000000000	126
0.675000000000000	87
0.725000000000000	127
0.775000000000000	109
0.825000000000000	115
0.875000000000000	98
0.925000000000000	118
0.975000000000000	95
1.02500000000000	128
1.07500000000000	77
1.12500000000000	48
1.17500000000000	34
1.22500000000000	37
1.27500000000000	44
1.32500000000000	27
1.37500000000000	34
1.42500000000000	33
1.47500000000000	32
1.52500000000000	37
1.57500000000000	24
1.62500000000000	29
1.67500000000000	37
1.72500000000000	46
1.77500000000000	53
1.82500000000000	38
1.87500000000000	44
1.92500000000000	43
1.97500000000000	285
}{\recCountTestDue}

\begin{axis}[%
  width=6.5cm,height=4cm,scale only axis,
  xmin=-2,xmax=2,ymin=-1,
  legend columns=1, 
  legend cell align={left},
  legend style={
        fill=none,
  	legend pos=south west,
  	draw=none,
        font=\tiny,
  	%cells={anchor=east}
  }
]
\addplot[ycomb,mark=none,blue] table[x=step,y=variation] 
{\recCountTestDue};
% \addlegendentry{exact};

%\draw[dashed,help lines] (axis cs:3.5,0.01) -- (axis cs:3.5,.12) ;
%\draw[dashed,help lines] (axis cs:8.5,0.01) -- (axis cs:8.5,.12) ;
%\draw[stealth-stealth,help lines] (axis cs:0.5,0.02) -- node[pos=0.5,below,black]{\tiny $m=50$} (axis cs:3.5,.02) ;
%\draw[stealth-stealth,help lines] (axis cs:3.5,0.02) -- node[pos=0.5,below,black]{\tiny $m=71$} (axis cs:8.5,.02) ;
%\draw[stealth-stealth,help lines] (axis cs:8.5,0.02) -- node[pos=0.5,below,black]{\tiny $m=100$} (axis cs:20.5,.02) ;

\end{axis}
\end{tikzpicture}
%\end{document} 
    \caption{Left: the exact solution of Test 2 in red and the numerical one in blue computed with a third order $\CWENOZ$ reconstruction and $N=81$. Right: reconstructions count in the last time step.}
    \label{fig:solTestDue}
\end{figure}
\begin{table}
\begin{center}
\begin{tabular}{|c|c|c|c|c|c|c|}
\hline &\multicolumn{2}{ |c|}{$\WENO$} &\multicolumn{2}{|c|}{$\CWENO$} &\multicolumn{2}{|c|}{$\CWENOZ$} \\
\hline $N$ & $L^{1}$-err & $L^1$-ord & $L^{1}$-err & $L^1$-ord & $L^{1}$-err & $L^1$-ord\\ \hline \hline
$81$ & $3.56e-06$ & & $2.24e-06$ & & $1.78e-06$ & \\ \hline 
$161$ & $2.83e-07$ & $3.65$ & $1.80e-07$ & $3.64$ & $1.44e-07$ & $3.63$\\ \hline 
$321$ & $2.45e-08$ & $3.53$ & $1.59e-08$ & $3.50$ & $1.30e-08$ & $3.47$\\ \hline 
$641$ & $2.61e-09$ & $3.23$ & $1.96e-09$ & $3.02$ & $1.75e-09$ & $2.90$\\ \hline 
% $1281$ & $7.27e-10$ & $1.84$ & $6.78e-10$ & $1.53$ & $6.63e-10$ & $1.40$\\ \hline 
% $2561$ & $5.69e-10$ & $0.35$ & $5.67e-10$ & $0.26$ & $5.67e-10$ & $0.23$\\ \hline 
% $5121$ & $5.66e-10$ & $0.01$ & $5.66e-10$ & $0.00$ & $5.65e-10$ & $0.00$\\ \hline 
% $10241$ & $5.45e-10$ & $0.06$ & $5.45e-10$ & $0.05$ & $5.45e-10$ & $0.05$\\ \hline 
  \end{tabular}
\end{center}
\caption{ Errors at time $T=1$ for Test 2, $\WENO$, $\CWENO$ and $\CWENOZ$  schemes.\label{tab:test2err} }
\end{table}

% CINECA
\begin{table}
\begin{center}
\begin{tabular}{|c|c|c|c|c|c|}
\hline & $\WENO$ & \multicolumn{2}{ |c|}{$\CWENO$} & \multicolumn{2}{ |c|}{$\CWENOZ$}\\
\hline grid & CPU time & CPU time & $\%$ gain & CPU time & $\%$ gain \\ \hline
\hline
$81$ & $6.19e-03$ & $4.83e-03$ & $21.95$ & $4.86e-03$ & $21.48$\\ \hline 
$161$ & $9.79e-03$ & $8.88e-03$ & $9.31$ & $8.64e-03$ & $11.76$\\ \hline 
$321$ & $2.74e-02$ & $2.39e-02$ & $13.03$ & $2.39e-02$ & $12.84$\\ \hline 
$641$ & $9.42e-02$ & $7.65e-02$ & $18.81$ & $7.68e-02$ & $18.41$\\ \hline 
$1281$ & $3.48e-01$ & $2.83e-01$ & $18.65$ & $2.90e-01$ & $16.72$\\ \hline 
$2561$ & $1.33e+00$ & $1.10e+00$ & $17.15$ & $1.10e+00$ & $16.83$\\ \hline 
$5121$ & $5.28e+00$ & $4.30e+00$ & $18.51$ & $4.29e+00$ & $18.77$\\ \hline 
$10241$ & $2.09e+01$ & $1.68e+01$ & $19.42$ & $1.68e+01$ & $19.46$\\ \hline 
  \end{tabular}
\end{center}

\caption{CPU times for Test 2, $\WENO$ and $\CWENO$ schemes.\label{tab:test2cpu} }
\end{table}

This test deals with the HJ equation
\begin{equation}\label{eq:test2}
\begin{cases}v_t(t,x)+\frac12|\partial _x v(t,x)|^2=0 \\
v(0,x)=v_0(x)=\min(-\cos(\pi x/2) ,0), 
\end{cases}
\end{equation}
in the domain $[-2,2]$, with homogeneous boundary conditions, which are treated using the extrapolation technique. The approximate solution is computed, at $T=1$, with $\ddt=10 \ddx$, while the exact solution $v(t,x)$ is defined as follows:
\begin{equation}
    v(t,x) = \min(0,\frac{1}{2}[a^*(t,x)]^2 + v_0(x-t a^*(t,x))),
\end{equation}
where, for a given $(t,x)$, $a^*(t,x)$  represents the optimal control, which  is constant along the characteristics, and it can be computed as the solution of 
\begin{equation}\label{eq:beta}
    a^*(t,x) = g(t,x,a^*(t,x))
\end{equation}
with the function $g$ defined by
\begin{equation}
    g(t,x,a^*(t,x)) = \begin{cases}
        -\frac{2}{\pi t}\arcsin{\frac{2a^*(t,x)}{\pi}} + \frac{x}{t} \quad &\text{if} \;|a^*(t,x)|\leq \frac{\pi}{2},\\
        \frac{\pi}{2} \quad &\text{if} \;a^*(t,x)> \frac{\pi}{2},\\
        -\frac{\pi}{2} \quad &\text{if} \;a^*(t,x)< -\frac{\pi}{2}.\\
    \end{cases}
\end{equation}
For our purposes, since the exact solution of \eqref{eq:beta} is not available, we compute  an approximate solution by using a fixed point algorithm with initial guess $a^*_0(t,x) = \frac{1}{2}x^2 t$.

Since the characteristics are affine,  in this test we have combined the third order spatial reconstructions with RK1 time stepping, which computes these curves exactly. Plots of the numerical and exact solutions are shown in the left panel of Figure \ref{fig:solTestDue} and $L^{1}$ errors are listed in Table \ref{tab:test2err}. 
In this case too, all schemes compute approximations affected by errors of the same order but the central reconstructions turn out to be more accurate (errors are about $2/3$ than the $\WENO$ case). In fact, the feet of the characteristics fall in regular regions of $v$ at each time step, and therefore all the reconstructions are close to the optimal cubic one. Table \ref{tab:test2cpu} reports the computational costs of each scheme, showing an average gain of $16-19\%$ for both the Central $\WENO$ reconstructions. This is again the effect of the multiple reconstructions performed in each cell (see the right panel of Figure~\ref{fig:solTestDue}, in which the number of reconstructions performed in each cell during the last time step is depicted). We point out that, as in Test 1, the high number of reconstructions computed in the first and in the last cell is related to the boundary conditions. 

\paragraph{Test 3: One--dimensional eikonal equation.}
\begin{figure}
\centering
%\documentclass{standalone}
%\usepackage{pgfplots}
%\begin{document}
\begin{tikzpicture}
\pgfplotstableread{%
step variation
0	-3.29309757410162e-05
0.0997331001139617	0.0995279771013144
0.199466200227923	0.198100608972365
0.299199300341885	0.294706562162566
0.398932400455847	0.388387001186603
0.498665500569808	0.478211905919010
0.598398600683770	0.563289314516671
0.698131700797732	0.642773251754288
0.797864800911694	0.715871886077055
0.897597901025655	0.781855619942927
0.997331001139617	0.840063991202133
1.09706410125358	0.889914318460071
1.19679720136754	0.930908857160448
1.29653030148150	0.962643867716495
1.39626340159546	0.984811842320849
1.49599650170943	0.997200853735506
1.59572960182339	0.999687506830722
1.69546270193735	0.992238198896621
1.79519580205131	0.974939219607758
1.89492890216527	0.947955268918204
1.99466200227923	0.911546093886080
2.09439510239320	0.866067123993358
2.19412820250716	0.811969971572238
2.29386130262112	0.749796068069264
2.39359440273508	0.680167831881691
2.49332750284904	0.603781351402214
2.59306060296300	0.521398092425894
2.69279370307697	0.433838125980593
2.79252680319093	0.341971392855681
2.89225990330489	0.246709914824939
2.99199300341885	0.148999326114348
3.09172610353281	0.0498094021016445
3.19145920364677	-0.0498752308868407
3.29119230376074	-0.149064824023784
3.39092540387470	-0.246774311545920
3.49065850398866	-0.342032522148198
3.59039160410262	-0.433890366012064
3.69012470421658	-0.521435742738000
3.78985780433054	-0.603799616513412
3.88959090444451	-0.680163676281667
3.98932400455847	-0.749769006822284
4.08905710467243	-0.811923861267286
4.18879020478639	-0.866010541484179
4.28852330490035	-0.911491519900672
4.38825640501432	-0.947914788339977
4.48798950512828	-0.974918360646801
4.58772260524224	-0.992233731592173
4.68745570535620	-0.999687619785267
4.78718880547016	-0.997200049837827
4.88692190558412	-0.984800051821619
4.98665500569808	-0.962613058896659
5.08638810581205	-0.930860161596835
5.18612120592601	-0.889857056373582
5.28585430603997	-0.840011263250655
5.38558740615393	-0.781818160791513
5.48532050626789	-0.715856075501488
5.58505360638186	-0.642780538994521
5.68478670649582	-0.563317794378240
5.78451980660978	-0.478257576975694
5.88425290672374	-0.388444761983173
5.98398600683770	-0.294769608070253
6.08371910695166	-0.198165691649591
6.18345220706563	-0.0995937046361612
}{\solTestTre}

\pgfplotstableread{%
step variation
0	0
0.0997331001139617	0.0995678465958167
0.199466200227923	0.198146143199398
0.299199300341885	0.294755174410904
0.398932400455847	0.388434796274695
0.498665500569808	0.478253978621318
0.598398600683770	0.563320058063622
0.698131700797732	0.642787609686539
0.797864800911694	0.715866849259718
0.897597901025655	0.781831482468030
0.997331001139617	0.840025923150771
1.09706410125358	0.889871808811469
1.19679720136754	0.930873748644204
1.29653030148150	0.962624246950012
1.39626340159546	0.984807753012208
1.49599650170943	0.997203797181180
1.59572960182339	0.999689182000816
1.69546270193735	0.992239206600172
1.79519580205131	0.974927912181824
1.89492890216527	0.947927346167132
1.99466200227923	0.911505852311673
2.09439510239320	0.866025403784439
2.19412820250716	0.811938005715856
2.29386130262112	0.749781202967734
2.39359440273508	0.680172737770919
2.49332750284904	0.603804410325477
2.59306060296300	0.521435203379498
2.69279370307697	0.433883739117558
2.79252680319093	0.342020143325669
2.89225990330489	0.246757397690293
2.99199300341885	0.149042266176174
3.09172610353281	0.0498458856606970
3.19145920364677	-0.0498458856606972
3.29119230376074	-0.149042266176174
3.39092540387470	-0.246757397690294
3.49065850398866	-0.342020143325669
3.59039160410262	-0.433883739117558
3.69012470421658	-0.521435203379498
3.78985780433054	-0.603804410325478
3.88959090444451	-0.680172737770920
3.98932400455847	-0.749781202967734
4.08905710467243	-0.811938005715857
4.18879020478639	-0.866025403784439
4.28852330490035	-0.911505852311673
4.38825640501432	-0.947927346167132
4.48798950512828	-0.974927912181824
4.58772260524224	-0.992239206600172
4.68745570535620	-0.999689182000816
4.78718880547016	-0.997203797181180
4.88692190558412	-0.984807753012208
4.98665500569808	-0.962624246950012
5.08638810581205	-0.930873748644204
5.18612120592601	-0.889871808811469
5.28585430603997	-0.840025923150771
5.38558740615393	-0.781831482468030
5.48532050626789	-0.715866849259718
5.58505360638186	-0.642787609686539
5.68478670649582	-0.563320058063622
5.78451980660978	-0.478253978621318
5.88425290672374	-0.388434796274695
5.98398600683770	-0.294755174410904
6.08371910695166	-0.198146143199398
6.18345220706563	-0.0995678465958164
}{\solTestTreExa}

\begin{axis}[%
  width=6.5cm,height=4cm,scale only axis,
  xmin=0,xmax=6.28318530717959,
  legend columns=1, 
  legend cell align={left},
  legend style={
        fill=none,
  	legend pos=south west,
  	draw=none,
        font=\tiny,
  	%cells={anchor=east}
  }
]
\addplot[color=red,mark=none] table[x=step,y=variation] {\solTestTreExa};
\addlegendentry{exact};
\addplot[color=blue,mark=star,only marks] table[x=step,y=variation] {\solTestTre};
\addlegendentry{numerical};

%\draw[dashed,help lines] (axis cs:3.5,0.01) -- (axis cs:3.5,.12) ;
%\draw[dashed,help lines] (axis cs:8.5,0.01) -- (axis cs:8.5,.12) ;
%\draw[stealth-stealth,help lines] (axis cs:0.5,0.02) -- node[pos=0.5,below,black]{\tiny $m=50$} (axis cs:3.5,.02) ;
%\draw[stealth-stealth,help lines] (axis cs:3.5,0.02) -- node[pos=0.5,below,black]{\tiny $m=71$} (axis cs:8.5,.02) ;
%\draw[stealth-stealth,help lines] (axis cs:8.5,0.02) -- node[pos=0.5,below,black]{\tiny $m=100$} (axis cs:20.5,.02) ;

\end{axis}
\end{tikzpicture}
%\end{document} 
%\documentclass{standalone}
%\usepackage{pgfplots}
%\begin{document}
\begin{tikzpicture}
\pgfplotstableread{%
step variation
0.0498665500569808	174
0.149599650170943	183
0.249332750284904	53
0.349065850398866	167
0.448798950512828	179
0.548532050626789	186
0.648265150740751	181
0.747998250854713	180
0.847731350968674	181
0.947464451082636	190
1.04719755119660	169
1.14693065131056	188
1.24666375142452	178
1.34639685153848	179
1.44612995165244	167
1.54586305176641	43
1.64559615188037	169
1.74532925199433	180
1.84506235210829	195
1.94479545222225	177
2.04452855233621	180
2.14426165245018	197
2.24399475256414	177
2.34372785267810	200
2.44346095279206	173
2.54319405290602	181
2.64292715301998	171
2.74266025313395	178
2.84239335324791	52
2.94212645336187	180
3.04185955347583	173
3.14159265358979	182
3.24132575370376	170
3.34105885381772	175
3.44079195393168	175
3.54052505404564	293
3.64025815415960	178
3.73999125427356	176
3.83972435438753	178
3.93945745450149	169
4.03919055461545	182
4.13892365472941	184
4.23865675484337	174
4.33838985495733	165
4.43812295507130	183
4.53785605518526	171
4.63758915529922	177
4.73732225541318	308
4.83705535552714	165
4.93678845564110	166
5.03652155575507	188
5.13625465586903	176
5.23598775598299	183
5.33572085609695	180
5.43545395621091	169
5.53518705632487	171
5.63492015643884	193
5.73465325655280	170
5.83438635666676	189
5.93411945678072	290
6.03385255689468	185
6.13358565700864	175
}{\recCountTestTre}

\begin{axis}[%
  width=6.5cm,height=4cm,scale only axis,
  xmin=0,xmax=6.28318530717959,
  legend columns=1, 
  legend cell align={left},
  legend style={
        fill=none,
  	legend pos=south west,
  	draw=none,
        font=\tiny,
  	%cells={anchor=east}
  }
]
\addplot[ycomb,mark=none,blue] table[x=step,y=variation] {\recCountTestTre};
% \addlegendentry{exact};

%\draw[dashed,help lines] (axis cs:3.5,0.01) -- (axis cs:3.5,.12) ;
%\draw[dashed,help lines] (axis cs:8.5,0.01) -- (axis cs:8.5,.12) ;
%\draw[stealth-stealth,help lines] (axis cs:0.5,0.02) -- node[pos=0.5,below,black]{\tiny $m=50$} (axis cs:3.5,.02) ;
%\draw[stealth-stealth,help lines] (axis cs:3.5,0.02) -- node[pos=0.5,below,black]{\tiny $m=71$} (axis cs:8.5,.02) ;
%\draw[stealth-stealth,help lines] (axis cs:8.5,0.02) -- node[pos=0.5,below,black]{\tiny $m=100$} (axis cs:20.5,.02) ;

\end{axis}
\end{tikzpicture}
%\end{document} 
    \caption{Left: the exact solution of Test 3 in red and the numerical one in blue computed with a third order $\CWENOZ$ scheme and $N=63$. Right: reconstructions count in the last time step.}
    \label{fig:solTestTre}
\end{figure}
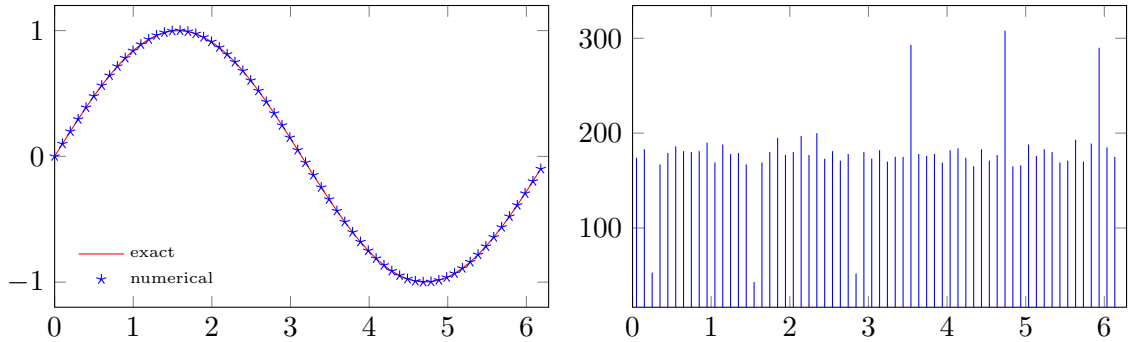
\begin{table}
\begin{center}
\begin{tabular}{|c|c|c|c|c|c|c|}
\hline &\multicolumn{2}{ |c|}{$\WENO$} &\multicolumn{2}{|c|}{$\CWENO$} &\multicolumn{2}{|c|}{$\CWENOZ$} \\
\hline $N$ & $L^{1}$-err & $L^1$-ord & $L^{1}$-err & $L^1$-ord & $L^{1}$-err & $L^1$-ord\\ \hline \hline
$126$ & $1.28e-05$ & & $1.38e-05$ & & $1.43e-05$ & \\ \hline 
$252$ & $1.54e-06$ & $3.05$ & $1.65e-06$ & $3.07$ & $1.70e-06$ & $3.07$\\ \hline 
$503$ & $1.93e-07$ & $3.00$ & $2.03e-07$ & $3.02$ & $2.10e-07$ & $3.02$\\ \hline 
$1006$ & $2.56e-08$ & $2.91$ & $2.63e-08$ & $2.95$ & $2.71e-08$ & $2.95$\\ \hline 
% $2011$ & $5.41e-09$ & $2.24$ & $5.41e-09$ & $2.28$ & $5.41e-09$ & $2.32$\\ \hline 
% $4022$ & $5.24e-09$ & $0.05$ & $4.92e-09$ & $0.14$ & $4.90e-09$ & $0.14$\\ \hline 
% $8043$ & $4.78e-09$ & $0.13$ & $4.78e-09$ & $0.04$ & $4.78e-09$ & $0.04$\\ \hline 
% $16085$ & $4.50e-09$ & $0.09$ & $4.48e-09$ & $0.09$ & $4.47e-09$ & $0.10$\\ \hline 
  \end{tabular}
\end{center}
\caption{ Errors at time $T=1$ for Test 3, $\WENO$, $\CWENO$ and $\CWENOZ$  schemes.\label{tab:test3err} }
\end{table}

% CINECA
\begin{table}
\begin{center}
\begin{tabular}{|c|c|c|c|c|c|}
\hline & $\WENO$ & \multicolumn{2}{ |c|}{$\CWENO$} & \multicolumn{2}{ |c|}{$\CWENOZ$}\\
\hline grid & CPU time & CPU time & $\%$ gain & CPU time & $\%$ gain \\ \hline
\hline
$126$ & $5.92e-02$ & $5.42e-02$ & $8.43$ & $5.64e-02$ & $4.81$\\ \hline 
$252$ & $2.11e-01$ & $1.94e-01$ & $8.04$ & $1.94e-01$ & $8.05$\\ \hline 
$503$ & $7.97e-01$ & $7.37e-01$ & $7.53$ & $7.48e-01$ & $6.05$\\ \hline 
$1006$ & $3.13e+00$ & $2.83e+00$ & $9.75$ & $2.87e+00$ & $8.45$\\ \hline 
$2011$ & $1.19e+01$ & $1.11e+01$ & $6.67$ & $1.11e+01$ & $6.79$\\ \hline 
$4022$ & $4.73e+01$ & $4.27e+01$ & $9.58$ & $4.21e+01$ & $11.03$\\ \hline 
$8043$ & $1.80e+02$ & $1.65e+02$ & $8.63$ & $1.62e+02$ & $10.02$\\ \hline 
$16085$ & $6.96e+02$ & $6.41e+02$ & $7.96$ & $6.38e+02$ & $8.35$\\ \hline 
  \end{tabular}
\end{center}
\caption{CPU times for Test 3, $\WENO$ and $\CWENO$ schemes.\label{tab:test3cpu} }
\end{table}

In this test, we consider the HJ equation
\begin{equation}\label{eq:test3}
\begin{cases}v_t(t,x)+\frac12|\partial _x v(t,x)|^2-f(t,x)=0 \\
v(0,x)=v_0(x), 
\end{cases}
\end{equation}
with $f(t,x)=-\sin(x)+(\frac{9}{8}+\frac{t^2-3t}{2})\cos^2(x)$ and $v_0(x)=\frac{3}{2}\sin(x)$.
We consider the problem in $[0,2\pi]$ with periodic boundary conditions and compute the numerical solution at time $T=0.5$, with $\ddt = \ddx$. This problem has a known exact solution given by $v(t,x)=(\frac{3}{2}-t)\sin(x)$.

In this case, characteristics are not affine and both the dynamic and the cost function explicitly depend on $x$ and $t$.
Therefore this is a specific test for the third order accuracy of our scheme, which depends on the use of the third order RK scheme \eqref{eq:RK3} and on the approximation of the integral of the cost function.

Plots of the exact and of the numerical solution are shown in the left panel of Fig.\ref{fig:solTestTre}, while in the right panel the reconstruction count in the last time step is depicted. We point out that in this case, with periodic boundary conditions, there is no artificial increase in the first and in the last cell near the boundary.

Errors and computational times are listed in Tables \ref{tab:test3err} and \ref{tab:test3cpu}. As in the previous tests, due to the regularity of the exact solution, all schemes perform almost equally when compared in terms of the error. Looking at the computational costs, the Central $\WENO$ reconstructions allow for a saving of about $6-10\%$.

\paragraph{Test 4: Two--dimensional semi-convex data.}
\begin{figure}
\centering
\epsfig{file=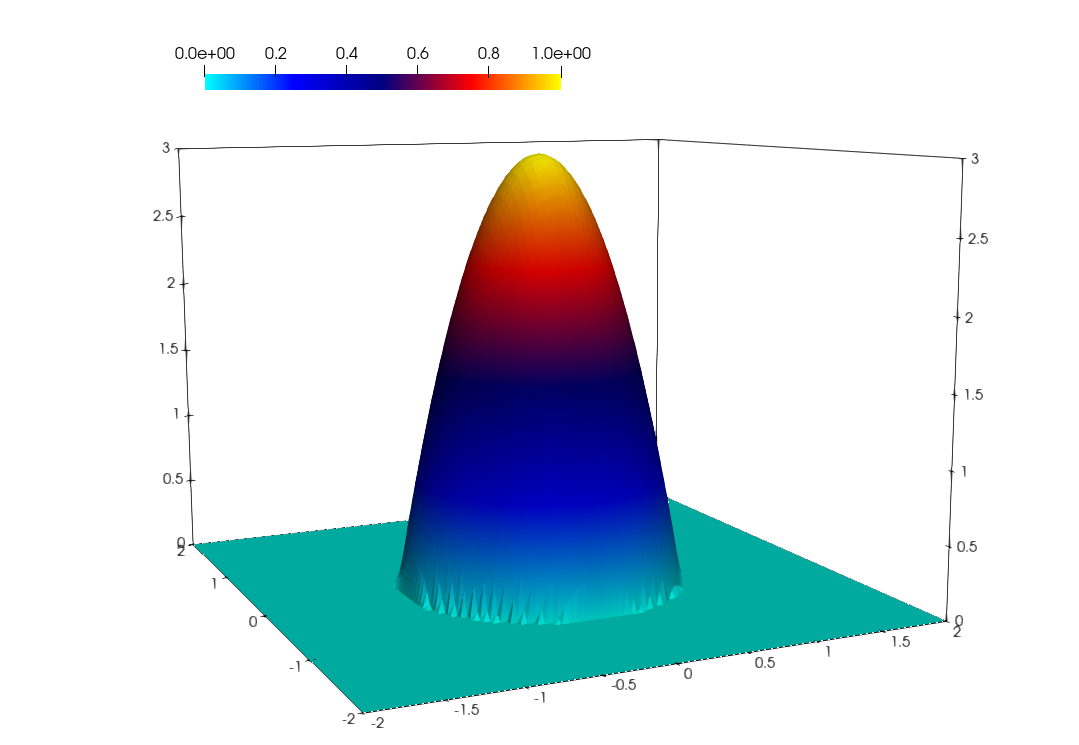,width=0.45\linewidth}\epsfig{file=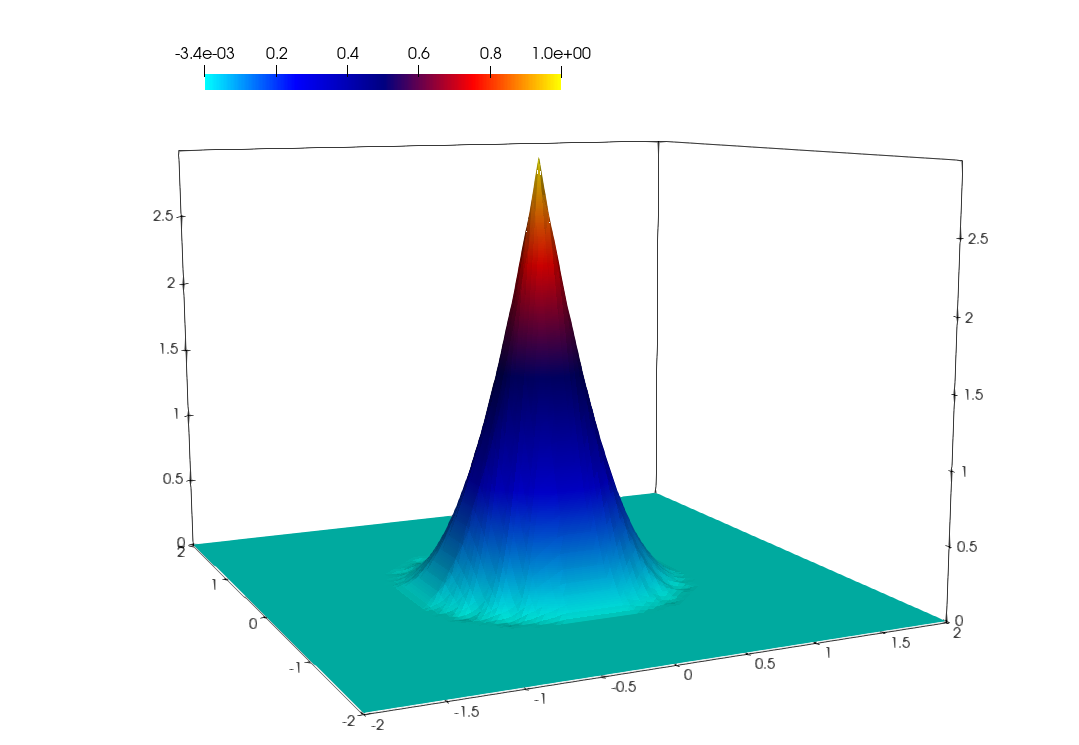,width=0.45\linewidth}
\caption{Test 4. Initial condition (left) and numerical solution computed with $\CWENO$ at $T=1$ (right) with $N=81$.}\label{fig:test4}
\end{figure}
% \begin{figure}
% \centering
% \begin{tabular}{cc}
% \input{grafici1d/test4radKinkCWENO}
% & \input{grafici1d/test4radKinkCWENOZ}
% \\
% \input{grafici1d/test4radMidCWENO}
% & \input{grafici1d/test4radMidCWENOZ}
% \\
% \input{grafici1d/test4radTailCWENO}
% & \input{grafici1d/test4radTailCWENOZ}
% \\
% \end{tabular}
%     \caption{Comparison of overshoots and undershoots of the numerical solutions computed with the traditional $\WENO$ scheme and with the central $\WENO$ ones on a $81\times81$ grid. In the left panels $\WENO$ and $\CWENO$ are compared, while in the right panels similar plots are shown for $\WENO$ and $\CWENOZ$. The exact solution is always represented with a thick black line. First row: scatter plots of all data points with $r\in[0,0.2]$ corresponding to the spike of the solution. Second row: scatter plots of all data points with $r\in[0.4,0.6]$ corresponding to regular region of the solution. Third row: scatter plots of all data points with $r\in[0.9,1.1]$ corresponding to the singular region of the initial data.}
%     \label{fig:scatterTestQuattro}
% \end{figure}

\begin{figure}
\centering
\begin{tabular}{cc}
%\documentclass{standalone}
%\usepackage{pgfplots}
%\begin{document}
\begin{tikzpicture}

\begin{axis}[%
  width=6cm,height=4cm,scale only axis,
  legend columns=1, 
  legend cell align={left},
  legend style={
        fill=none,
  	legend pos=north east,
  	draw=none,
        font=\tiny,
    enlarge x limits=0,
  	%cells={anchor=east}
  }
]
%\addplot[color=red,scatter] table[y=zz,x=rr]{datiRadialiWENO321.txt};
%\addlegendentry{$\WENO$};
\addplot[color=black,thick] table[col sep=comma]{grafici1d/exa81kink.txt};
\addlegendentry{exact};
\addplot[color=red,only marks,mark size=2pt,mark=x] table[col sep=comma]{grafici1d/cweno161kink.txt};
\addlegendentry{$\CWENO$};
\addplot[color=blue,only marks,mark size=0.5pt,mark=o] table[col sep=comma]{grafici1d/weno161kink.txt};
\addlegendentry{$\WENO$};

\end{axis}
\end{tikzpicture}
%\end{document} 
& %\documentclass{standalone}
%\usepackage{pgfplots}
%\begin{document}
\begin{tikzpicture}

\begin{axis}[%
  width=6cm,height=4cm,scale only axis,
  legend columns=1, 
  legend cell align={left},
  legend style={
        fill=none,
  	legend pos=north east,
  	draw=none,
        font=\tiny,
    enlarge x limits=0,
  	%cells={anchor=east}
  }
]
%\addplot[color=red,scatter] table[y=zz,x=rr]{datiRadialiWENO321.txt};
%\addlegendentry{$\WENO$};
\addplot[color=black,thick] table[col sep=comma]{grafici1d/exa81kink.txt};
\addlegendentry{exact};
\addplot[color=red,only marks,mark size=2pt,mark=x] table[col sep=comma]{grafici1d/cwenoz161kink.txt};
\addlegendentry{$\CWENO$};
\addplot[color=blue,only marks,mark size=0.5pt,mark=o] table[col sep=comma]{grafici1d/weno161kink.txt};
\addlegendentry{$\CWENOZ$};

\end{axis}
\end{tikzpicture}
%\end{document} 
\\
%\documentclass{standalone}
%\usepackage{pgfplots}
%\begin{document}
\begin{tikzpicture}

\begin{axis}[%
  width=6cm,height=4cm,scale only axis,
  legend columns=1, 
  legend cell align={left},
  legend style={
        fill=none,
  	legend pos=north east,
  	draw=none,
        font=\tiny,
    enlarge x limits=0,
  	%cells={anchor=east}
  }
]
%\addplot[color=red,scatter] table[y=zz,x=rr]{datiRadialiWENO321.txt};
%\addlegendentry{$\WENO$};
\addplot[color=black,thick] table[col sep=comma]{grafici1d/exa81mid.txt};
\addlegendentry{exact};
\addplot[color=red,only marks,mark size=2pt,mark=x] table[col sep=comma]{grafici1d/cweno81mid.txt};
\addlegendentry{$\CWENO$};
\addplot[color=blue,only marks,mark size=0.5pt,mark=o] table[col sep=comma]{grafici1d/weno81mid.txt};
\addlegendentry{$\WENO$};

\end{axis}
\end{tikzpicture}
%\end{document} 
& %\documentclass{standalone}
%\usepackage{pgfplots}
%\begin{document}
\begin{tikzpicture}

\begin{axis}[%
  width=6cm,height=4cm,scale only axis,
  legend columns=1, 
  legend cell align={left},
  legend style={
        fill=none,
  	legend pos=north east,
  	draw=none,
        font=\tiny,
    enlarge x limits=0,
  	%cells={anchor=east}
  }
]
%\addplot[color=red,scatter] table[y=zz,x=rr]{datiRadialiWENO321.txt};
%\addlegendentry{$\WENO$};
\addplot[color=black,thick] table[col sep=comma]{grafici1d/exa81mid.txt};
\addlegendentry{exact};
\addplot[color=red,only marks,mark size=2pt,mark=x] table[col sep=comma]{grafici1d/cwenoz81mid.txt};
\addlegendentry{$\CWENOZ$};
\addplot[color=blue,only marks,mark size=0.5pt,mark=o] table[col sep=comma]{grafici1d/weno81mid.txt};
\addlegendentry{$\WENO$};

\end{axis}
\end{tikzpicture}
%\end{document} 
\\
%\documentclass{standalone}
%\usepackage{pgfplots}
%\begin{document}
\begin{tikzpicture}

\begin{axis}[%
  width=6cm,height=4cm,scale only axis,
  legend columns=1, 
  legend cell align={left},
  legend style={
        fill=none,
  	legend pos=north east,
  	draw=none,
        font=\tiny,
    enlarge x limits=0,
  	%cells={anchor=east}
  }
]
%\addplot[color=red,scatter] table[y=zz,x=rr]{datiRadialiWENO321.txt};
%\addlegendentry{$\WENO$};
\addplot[color=black,thick] table[col sep=comma]{grafici1d/exa81tail.txt};
\addlegendentry{exact};
\addplot[color=red,only marks,mark size=2pt,mark=x] table[col sep=comma]{grafici1d/cweno161tail.txt};
\addlegendentry{$\CWENO$};
\addplot[color=blue,only marks,mark size=0.5pt,mark=o] table[col sep=comma]{grafici1d/weno161tail.txt};
\addlegendentry{$\WENO$};

\end{axis}
\end{tikzpicture}
%\end{document} 
& %\documentclass{standalone}
%\usepackage{pgfplots}
%\begin{document}
\begin{tikzpicture}

\begin{axis}[%
  width=6cm,height=4cm,scale only axis,
  legend columns=1, 
  legend cell align={left},
  legend style={
        fill=none,
  	legend pos=north east,
  	draw=none,
        font=\tiny,
    enlarge x limits=0,
  	%cells={anchor=east}
  }
]
%\addplot[color=red,scatter] table[y=zz,x=rr]{datiRadialiWENO321.txt};
%\addlegendentry{$\WENO$};
\addplot[color=black,thick] table[col sep=comma]{grafici1d/exa81tail.txt};
\addlegendentry{exact};
\addplot[color=red,only marks,mark size=2pt,mark=x] table[col sep=comma]{grafici1d/cwenoz161tail.txt};
\addlegendentry{$\CWENOZ$};
\addplot[color=blue,only marks,mark size=0.5pt,mark=o] table[col sep=comma]{grafici1d/weno161tail.txt};
\addlegendentry{$\WENO$};

\end{axis}
\end{tikzpicture}
%\end{document} 
\\
\end{tabular}
    \caption{Test 4. Comparison of overshoots and undershoots of the numerical solutions computed with the traditional $\WENO$ scheme and with the central $\WENO$ ones on a $161\times161$ grid. In the left panels $\WENO$ and $\CWENO$ are compared, while in the right panels similar plots are shown for $\WENO$ and $\CWENOZ$. The exact solution is always represented with a thick black line. First row: scatter plots of all data points with $r\in[0,0.2]$ corresponding to the spike of the solution. Second row: scatter plots of all data points with $r\in[0.4,0.6]$ corresponding to regular region of the solution. Third row: scatter plots of all data points with $r\in[0.9,1.1]$ corresponding to the singular region of the initial data.}
    \label{fig:scatterTestQuattro}
\end{figure}
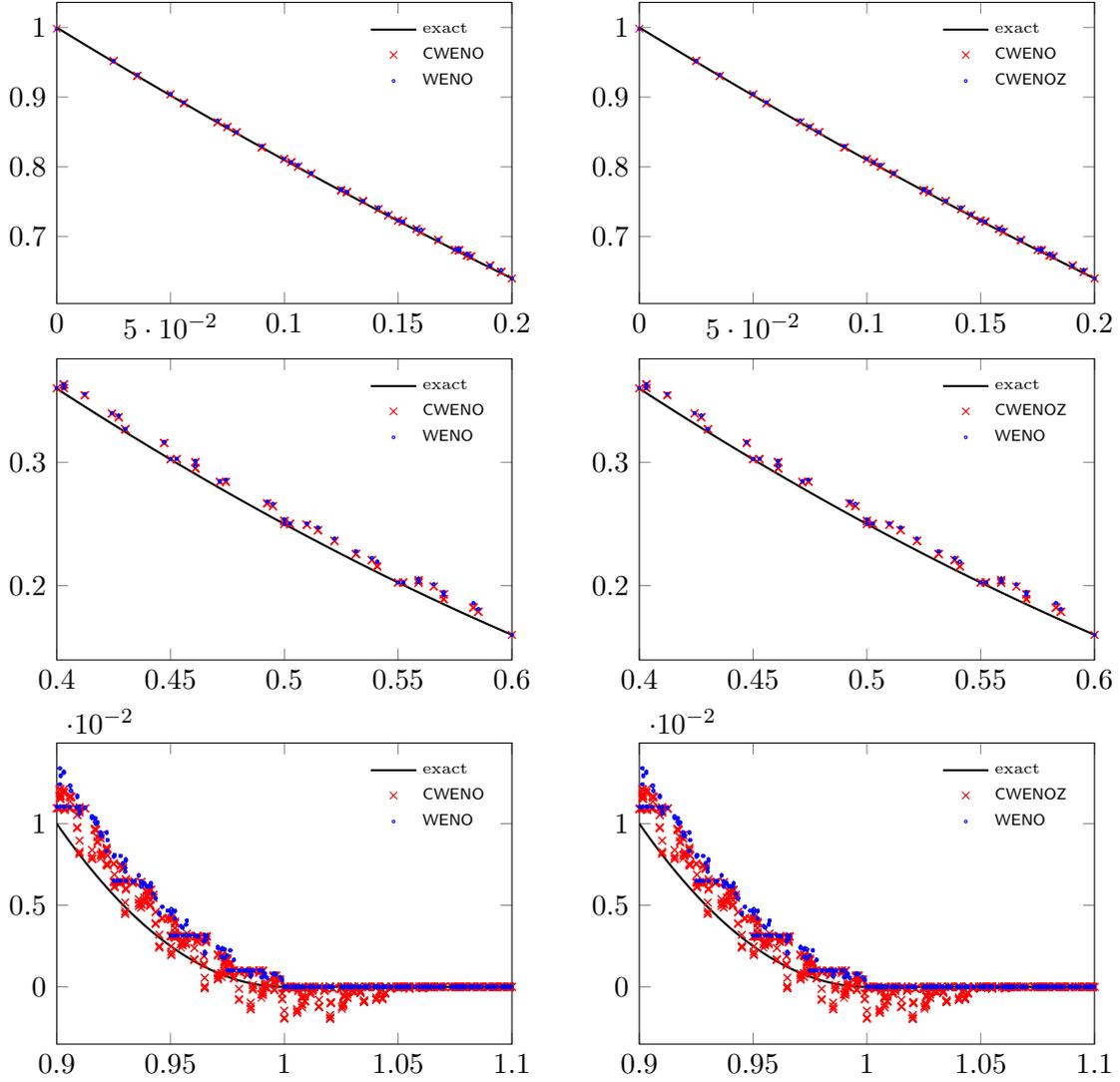

\begin{table}
\begin{center}
\begin{tabular}{|c|c|c|c|c|c|c|}
\hline &\multicolumn{2}{ |c|}{$\WENO$} &\multicolumn{2}{|c|}{$\CWENO$} &\multicolumn{2}{|c|}{$\CWENOZ$} \\
\hline $N$ & $L^{1}$-err & $L^1$-ord & $L^{1}$-err & $L^1$-ord & $L^{1}$-err & $L^1$-ord\\ \hline \hline
$41\times41$ & $4.02e-02$ & & $3.38e-02$ & & $3.38e-02$ & \\ \hline 
$81\times81$ & $2.25e-02$ & $0.84$ & $1.82e-02$ & $0.90$ & $1.81e-02$ & $0.90$\\ \hline 
$161\times161$ & $1.15e-02$ & $0.96$ & $9.01e-03$ & $1.01$ & $8.99e-03$ & $1.01$\\ \hline 
$321\times321$ & $5.56e-03$ & $1.05$ & $4.10e-03$ & $1.14$ & $4.09e-03$ & $1.14$\\ \hline 
$641\times641$ & $2.85e-03$ & $0.97$ & $2.03e-03$ & $1.02$ & $2.02e-03$ & $1.02$\\ \hline 
  \end{tabular}
\end{center}
\caption{ Errors at time $T=1$ for Test 4, $\WENO$, $\CWENO$ and $\CWENOZ$  schemes.\label{tab:test4err} }
\end{table}

\begin{table}
\begin{center}
\begin{tabular}{|c|c|c|c|c|c|}
\hline & $\WENO$ & \multicolumn{2}{ |c|}{$\CWENO$} & \multicolumn{2}{ |c|}{$\CWENOZ$}\\
\hline grid & CPU time & CPU time & $\%$ gain & CPU time & $\%$ gain \\ \hline
\hline
$41\times41$ & $2.88e+00$ & $1.98e+00$ & $31.37$ & $1.96e+00$ & $32.07$\\ \hline 
$81\times81$ & $2.19e+01$ & $1.54e+01$ & $29.70$ & $1.53e+01$ & $30.07$\\ \hline 
$161\times161$ & $1.72e+02$ & $1.19e+02$ & $31.10$ & $1.18e+02$ & $31.35$\\ \hline 
$321\times321$ & $1.37e+03$ & $9.38e+02$ & $31.56$ & $9.34e+02$ & $31.86$\\ \hline 
$641\times641$ & $1.09e+04$ & $7.44e+03$ & $31.59$ & $7.51e+03$ & $30.93$\\ \hline 
\end{tabular}
\end{center}
\caption{CPU times for Test 4, $\WENO$ and $\CWENO$ schemes.\label{tab:test4cpu}}
\end{table}

\begin{table}
\begin{center}
\small
\begin{tabular}{|c|c|c|c|c|c|c|}
\hline &\multicolumn{2}{ |c|}{$\WENO$} &\multicolumn{2}{|c|}{$\CWENO$} &\multicolumn{2}{|c|}{$\CWENOZ$} \\
\hline $N$ & min & max & min & max & min & max\\ \hline \hline
$21\times21$ & $-5.67e-05$ & $9.76e-01$ & $-6.08e-03$ & $9.80e-01$ & $-6.08e-03$ & $9.80e-01$\\ \hline 
$41\times41$ & $-4.62e-05$ & $9.89e-01$ & $-5.19e-03$ & $9.92e-01$ & $-5.19e-03$ & $9.92e-01$\\ \hline 
$81\times81$ & $-1.21e-05$ & $9.95e-01$ & $-3.39e-03$ & $9.96e-01$ & $-3.39e-03$ & $9.96e-01$\\ \hline 
$161\times161$ & $-3.45e-06$ & $9.98e-01$ & $-1.97e-03$ & $9.98e-01$ & $-1.97e-03$ & $9.98e-01$\\ \hline 
  \end{tabular}
\end{center}
\caption{Overshoots and undershoots observed in Test 4, $\WENO$, $\CWENO$ and $\CWENOZ$  schemes.\label{tab:test4undershoot} }
\end{table}

In this test, the HJ equation
\begin{equation}\label{eq:test4}
\begin{cases}v_t(t,x)+\frac12|Dv(t,x)|^2=0 \\
v(0,x)=v_0(x)=\max(1-|x|^2,0), 
\end{cases}
\end{equation}
is considered in $[-2,2]^2$, using extrapolation techniques to treat the homogeneous boundary condition. The exact solution of this problem is known and is given, for $t\geq1/2$, by
\begin{equation}
    v(t,x) = \begin{cases}
        \frac{(|x|-1)^2}{2t} \quad &\text{if} \; |x|\leq 1,\\
        0 \quad &\text{if} \; |x|>1.
    \end{cases}
\end{equation}
We set final time $T=0.5$ and compute the approximate solution with $\ddt=\frac{5}{4}\ddx$. Initial condition and final numerical solution are shown in Figure \ref{fig:test4}, 
while $L^1$ errors and CPU times are reported respectively in Tables \ref{tab:test4err} and \ref{tab:test4cpu}.

Since in this case the feet of characteristics always fall in the singular region of the solution, we observe, as expected, that all numerical schemes are degraded to first order accuracy. From Table~\ref{tab:test4err} one can see that the $L^1$-norm of the error is about $30\%$ lower for the central reconstructions, as in the previous cases. From Table~\ref{tab:test4undershoot} and Figure~\ref{fig:scatterTestQuattro} we can see that near the kink in the solution, the central reconstructions produce larger undershoots in the flat region, while being closer to the exact solution in the nonzero region. Finally, from Table \ref{tab:test4cpu} one can observe that central reconstructions schemes are about $30\%$ faster.

\paragraph{Test 5:  Front propagation with obstacles.}
\begin{figure}
\centering
\epsfig{file=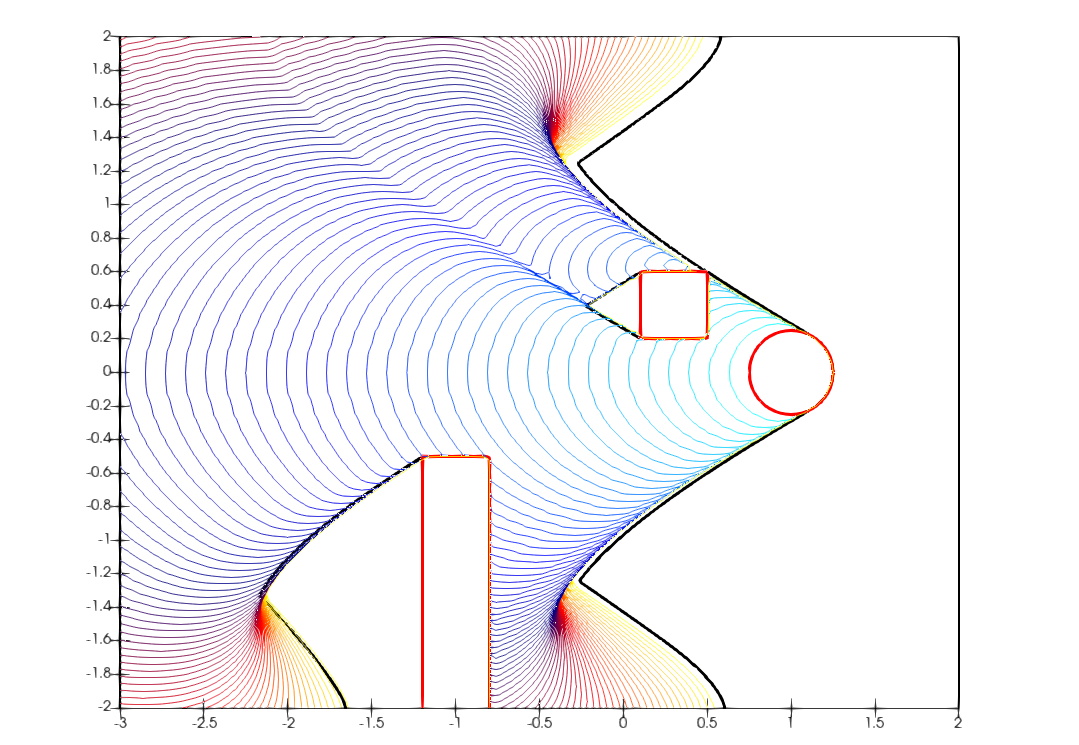,width=8cm}\epsfig{file=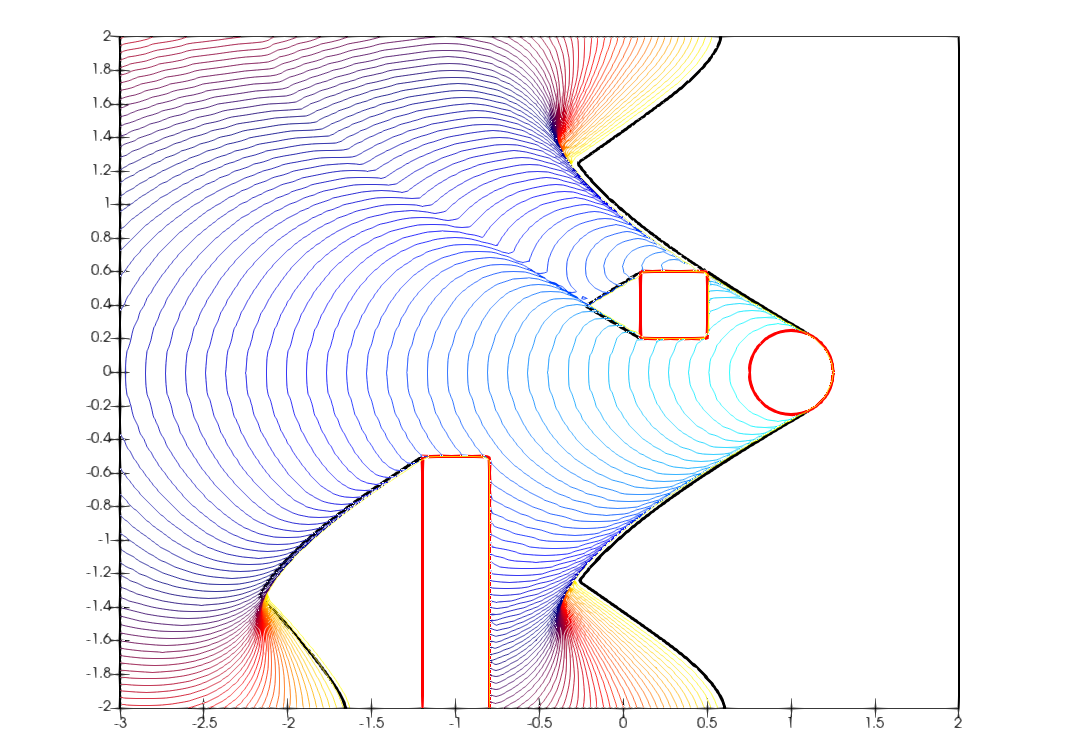,width=8cm}
\caption{Reachable sets \eqref{def:Rset} computed at each time step for Test 5 on a grid $101\times 81$ and $T=3$, using $\WENO$ (left) and $\CWENO$ (right) reconstructions. In both panel the black line represents a reference solution computed on a grid $1001\times 801$, the red circle represents the target and the red rectangles represent the obstacles.}\label{fig:test5}
\end{figure}

% CINECA
\begin{table}
\begin{center}
\begin{tabular}{|c|c|c|c|c|c|}
\hline & $\WENO$ & \multicolumn{2}{ |c|}{$\CWENO$} & \multicolumn{2}{ |c|}{$\CWENOZ$}\\
\hline grid & CPU time & CPU time & $\%$ gain & CPU time & $\%$ gain \\ \hline
\hline
$126\times126$ & $6.28e+00$ & $5.34e+00$ & $14.94$ & $5.30e+00$ & $15.66$\\ \hline 
$251\times251$ & $4.56e+01$ & $3.96e+01$ & $13.04$ & $3.96e+01$ & $13.06$\\ \hline 
$501\times501$ & $3.61e+02$ & $3.11e+02$ & $13.87$ & $3.20e+02$ & $11.50$\\ \hline 
$1001\times1001$ & $2.81e+03$ & $2.43e+03$ & $13.64$ & $2.40e+03$ & $14.59$\\ \hline 
$2001\times2001$ & $2.20e+04$ & $1.92e+04$ & $13.01$ & $1.89e+04$ & $14.00$\\ \hline 
  \end{tabular}
\end{center}
\caption{CPU times for Test 5, $\WENO$ and $\CWENO$ schemes.\label{tab:test5cpu} }
\end{table}

Last, we consider a state constrained Zermelo problem,
proposed in \cite{BFZ10}, for a swimmer that has to reach a circular island at $(1,0)$, facing a non-constant current flowing in the horizontal direction. Two obstacles are placed in the domain, near the target.

More precisely, to consider the minimum time problem the running cost $f_C(t,x,a)$ is set to zero, and the nonlinear dynamics is given by
$f_D(t,x,a)=(2-0.5x_2^2+a_1,a_2)$, with $A:=\{a \in \mathbb{R}^2 \,|\, \|a\|=(a_1^2+a^2_2)^{1/2}= 1\}$.
The swimmer has unit speed and can move in any direction $a\in A$; the current has speed $2-0.5x_2^2$ and flows in the direction $(1,0)$.
The target is the set $\mathcal{X}=\{x \;|\; v_0(x)\leq 0\}$, where $v_0(x)=C\,{\rm{min}}\big(\|(x_1-1,x_2)\|-r,r \big)$ with $r=0.25$ and $C=20$.
In order to take into account the presence of the two obstacles, we follow the level set approach, as proposed in \cite{BFZ10}, and represent the set of constrained states as
$\mathcal{K}=\{x \in \mathbb{R}^2\; | \; g(x)\leq 0\}$ where
$$ g(x)={\rm{max}}\Big(-\gamma, C\big( \gamma-{\rm{max}} (|x_1-0.3|,|x_2-0.4|)\big), C\big( \gamma-{\rm{max}}(|x_1+1|,|x_2+1.5|) \big ) \Big),$$
with $C=20$, $\gamma=0.2$.

For every $t\in[0,T]$, the set of points from which the swimmer can reach the
target {\em before time $t$} (Backward Reachable set) is represented by the level set 
$\mathcal{R}[0,t]=\{x\in \mathbb{R}^2\,|\, \exists\, s \in [0,t],\, v(x,s)\leq 0\}$ where $v$ is the solution to a constrained version of \eqref{eq:HJ} 
\begin{equation}\label{eq:HJC}
 \min(v_t+H(t,x,Dv),v-g(x))=0
\end{equation} 
% \todo{serve veramente? Non basta dire qui che $f_C=0$?}
% %\underset{a \in \mahtbb{R}^2, \|a\| \leq 1, b\in[0,1]}
% $$H(x,D v):=\underset{a \in A}
% {\rm{max}} \Big(-f_D(x,a)\cdot Dv \Big)=
% {\rm{max}}\Big(0, -(c-0.5x_2^2)\partial_{x_1}v +|Dv| \Big),
% $$
with initial condition $v(x,0)={\rm{max}}(v_0(x),g(x))$, (see \cite[Remark~2]{BFZ10}).
The space domain is defined as $[-3,2]\times[-2,2]$ and the final time to $T=3$.
By using scheme \eqref{eq:SL}, we discetize \eqref{eq:HJC} as
$$ \min(u^{n+1}_i-\underset{\underline{a} \in A^\nu}\min \{R[u^n](y^{n}_i(\underline{a}))\} ,u^{n+1}-g(x_i))=0,
$$
resulting in the following time marching scheme
\[
{u}_i^{n+1} = \max\left( \underset{\underline{a} \in A^\nu}\min \{R[u^n](y^{n}_i(\underline{a}))\} , g(x_i)\right),
\]
from which we define the numerical approximation of the seachable set $\mathcal{R}[0,t_n]$ as
%We observe that in this problem the  boundary is outflow, and the scheme auto  boundary conditions 
%we compute the numerical solution at time $t^{n+1}$ for node $x_i$ as
%
\begin{equation}\label{def:Rset}
    \mathcal{R}^n=\{x_i\,|\, \exists\, k \in \{0,\dots,n\}, u^k_i\leq 0\}.
\end{equation}

%where $\widetilde{v}_i^{n+1}$ is the solution of the numerical scheme \eqref{eq:SL} for \eqref{eq:HJ}.

The solution is computed with RK3 and $\ddt= \ddx$.
In Figure~\ref{fig:test5}, we show that, even in this very complex situation, no relevant difference in the accuracy of the computed solution can be observed between the $\WENO$ and the Central $\WENO$ approaches. On the other hand, Central reconstruction schemes provide a reduction of about $13-14\%$ of the computational time.

\section*{Conclusion}

In this study, a high order semi-Lagrangian scheme has been developed to approximate first order Hamilton--Jacobi--Bellman equations. The primary challenge lies in the nonsmooth nature of the solutions, which would lead to spurious oscillations when using unlimited high order polynomial interpolation.

To address this issue, the reconstruction of the solution at the previous time step at the foot of characteristics has to be performed with a non-oscillatory scheme. Since one has to perform very many interpolations during the minimization procedure, a $\CWENO$ technique was proposed to be coupled with the semi-Lagrangian approach.
Our study demonstrates that, in terms of errors, our new scheme maintains the favorable behavior of $\WENO$ schemes of \cite{CFR05}, producing about $30\%$ more accurate results in smooth regions at the price of some extra over/undershoots. However, its computational cost is significantly lower, by a $10-30\%$ depending on the specific test.

Additionally, we have established a convergence result in a simpler case, and provided several numerical simulations to further validate the effectiveness of our proposed scheme.

\subsection*{Future Directions}

Future research directions could explore the inclusion of high order treatments for boundary conditions and extensions of the convergence analysis to more general Hamiltonians. Investigating the scheme's performance in more complex scenarios would also contribute to a comprehensive understanding of its capabilities.

\subsection*{Acknowledgments}
The first two authors would like to thank the Italian Ministry of Instruction, University and Research (MIUR) for supporting this research with funds coming from
the PRIN Project 2022 (2022238YY5, entitled "Optimal control problems: analysis,
approximation").
All authors are members of the GNCS-INdAM Gruppo Nazionale per il Calcolo Scientifico of Istituto Nazionale di Alta Matematica, Francesco Severi, P.le Aldo Moro, Roma, Italy.

\bibliographystyle{alpha}
\bibliography{cweno_x}
\section*{Appendix}\label{sec:appendix}
\setcounter{MaxMatrixCols}{16}
We report here the matrix involved in the definition \eqref{eq:ind2dCoeff}.

\rotatebox{90}
{
\renewcommand*{\arraystretch}{1.5}
$
M = \begin{pmatrix} 
0 & 0 & 0 & 0 & 0 & 0 & 0 & 0 & 0 & 0 & 0 & 0 & 0 & 0 & 0 & 0 \\
0 & 0 & 0 & 0 & 0 & 0 & 0 & 0 & 0 & 0 & 0 & 0 & 0 & 0 & 0 & 0 \\
0 & 0 & 0 & 0 & 0 & 0 & 0 & 0 & 0 & 0 & 0 & 0 & 0 & 0 & 0 & 0 \\
0 & 0 & 0 & 4 & 0 & 0 & 6 & 2 & 0 & 0 & 3 & \frac{4}{3} & 0 & 2 & 1 & \frac{3}{2} \\
0 & 0 & 0 & 0 & 1 & 0 & 0 & 1 & 1 & 0 & 1 & 1 & 1 & 1 & 1 & 1 \\
0 & 0 & 0 & 0 & 0 & 4 & 0 & 0 & 2 & 6 & 0 & \frac{4}{3} & 3 & 1 & 2 & \frac{3}{2} \\
0 & 0 & 0 & 6 & 0 & 0 & 84 & 3 & 0 & 0 & 42 & 2 & 0 & 28 & \frac{3}{2} & 21 \\
0 & 0 & 0 & 2 & 1 & 0 & 3 & \frac{32}{3} & 1 & 0 & \frac{31}{2} & \frac{31}{3} & 1 & 15 & \frac{152}{15} & \frac{147}{10} \\
0 & 0 & 0 & 0 & 1 & 2 & 0 & 1 & \frac{32}{3} & 3 & 1 & \frac{31}{3} & \frac{31}{2} & \frac{152}{15} & 15 & \frac{147}{10} \\
0 & 0 & 0 & 0 & 0 & 6 & 0 & 0 & 3 & 84 & 0 & 2 & 42 & \frac{3}{2} & 28 & 21 \\
0 & 0 & 0 & 3 & 1 & 0 & 42 & \frac{31}{2} & 1 & 0 & \frac{989}{5} & 15 & 1 & \frac{954}{5} & \frac{147}{10} & \frac{933}{5} \\
0 & 0 & 0 & \frac{4}{3} & 1 & \frac{4}{3} & 2 & \frac{31}{3} & \frac{31}{3} & 2 & 15 & \frac{3992}{45} & 15 & \frac{1918}{15} & \frac{1918}{15} & \frac{737}{4} \\
0 & 0 & 0 & 0 & 1 & 3 & 0 & 1 & \frac{31}{2} & 42 & 1 & 15 & \frac{989}{5} & \frac{147}{10} & \frac{954}{5} & \frac{933}{5} \\
0 & 0 & 0 & 2 & 1 & 1 & 28 & 15 & \frac{152}{15} & \frac{3}{2} & \frac{954}{5} & \frac{1918}{15} & \frac{147}{10} & \frac{56076}{35} & \frac{737}{4} & \frac{161501}{70} \\
0 & 0 & 0 & 1 & 1 & 2 & \frac{3}{2} & \frac{152}{15} & 15 & 28 & \frac{147}{10} & \frac{1918}{15} & \frac{954}{5} & \frac{737}{4} & \frac{56076}{35} & \frac{161501}{70} \\
0 & 0 & 0 & \frac{3}{2} & 1 & \frac{3}{2} & 21 & \frac{147}{10} & \frac{147}{10} & 21 & \frac{933}{5} & \frac{737}{4} & \frac{933}{5} & \frac{161501}{70} & \frac{161501}{70} & \frac{721401}{25} \\
\end{pmatrix}.
$
} 

Also, we report the matrices involved in the definitions \eqref{eq:ind2dDati}. We assume that the vector $U$ of data in the $4\times4$ stencil is ordered lexicographically, i.e. $x$ direction faster and $y$ direction slower.

\rotatebox{90}
{\renewcommand*{\arraystretch}{1.5}
$
A_{\text{ne}} = \begin{pmatrix}
0 & 0 & 0 & 0 & 0 & 0 & 0 & 0 & 0 & 0 & 0 & 0 & 0 & 0 & 0 & 0 \\
0 & 0 & 0 & 0 & 0 & 0 & 0 & 0 & 0 & 0 & 0 & 0 & 0 & 0 & 0 & 0 \\
0 & 0 & 0 & 0 & 0 & 0 & 0 & 0 & 0 & 0 & 0 & 0 & 0 & 0 & 0 & 0 \\
0 & 0 & 0 & 0 & 0 & 0 & 0 & 0 & 0 & 0 & 0 & 0 & 0 & 0 & 0 & 0 \\
0 & 0 & 0 & 0 & 0 & 0 & 0 & 0 & 0 & 0 & 0 & 0 & 0 & 0 & 0 & 0 \\
0 & 0 & 0 & 0 & 0 & \frac{3497}{720} & \frac{-623}{90} & \frac{1787}{720} & 0 & \frac{-623}{90} & \frac{1667}{180} & \frac{-571}{180} & 0 & \frac{1787}{720} & \frac{-571}{180} & \frac{797}{720} \\
0 & 0 & 0 & 0 & 0 & \frac{-623}{90} & \frac{1141}{90} & \frac{-229}{45} & 0 & \frac{1667}{180} & \frac{-1547}{90} & \frac{1187}{180} & 0 & \frac{-571}{180} & \frac{278}{45} & \frac{-421}{180} \\
0 & 0 & 0 & 0 & 0 & \frac{1787}{720} & \frac{-229}{45} & \frac{1817}{720} & 0 & \frac{-571}{180} & \frac{1187}{180} & \frac{-293}{90} & 0 & \frac{797}{720} & \frac{-421}{180} & \frac{827}{720} \\
0 & 0 & 0 & 0 & 0 & 0 & 0 & 0 & 0 & 0 & 0 & 0 & 0 & 0 & 0 & 0 \\
0 & 0 & 0 & 0 & 0 & \frac{-623}{90} & \frac{1667}{180} & \frac{-571}{180} & 0 & \frac{1141}{90} & \frac{-1547}{90} & \frac{278}{45} & 0 & \frac{-229}{45} & \frac{1187}{180} & \frac{-421}{180} \\
0 & 0 & 0 & 0 & 0 & \frac{1667}{180} & \frac{-1547}{90} & \frac{1187}{180} & 0 & \frac{-1547}{90} & \frac{1472}{45} & \frac{-1157}{90} & 0 & \frac{1187}{180} & \frac{-1157}{90} & \frac{887}{180} \\
0 & 0 & 0 & 0 & 0 & \frac{-571}{180} & \frac{1187}{180} & \frac{-293}{90} & 0 & \frac{278}{45} & \frac{-1157}{90} & \frac{571}{90} & 0 & \frac{-421}{180} & \frac{887}{180} & \frac{-109}{45} \\
0 & 0 & 0 & 0 & 0 & 0 & 0 & 0 & 0 & 0 & 0 & 0 & 0 & 0 & 0 & 0 \\
0 & 0 & 0 & 0 & 0 & \frac{1787}{720} & \frac{-571}{180} & \frac{797}{720} & 0 & \frac{-229}{45} & \frac{1187}{180} & \frac{-421}{180} & 0 & \frac{1817}{720} & \frac{-293}{90} & \frac{827}{720} \\
0 & 0 & 0 & 0 & 0 & \frac{-571}{180} & \frac{278}{45} & \frac{-421}{180} & 0 & \frac{1187}{180} & \frac{-1157}{90} & \frac{887}{180} & 0 & \frac{-293}{90} & \frac{571}{90} & \frac{-109}{45} \\
0 & 0 & 0 & 0 & 0 & \frac{797}{720} & \frac{-421}{180} & \frac{827}{720} & 0 & \frac{-421}{180} & \frac{887}{180} & \frac{-109}{45} & 0 & \frac{827}{720} & \frac{-109}{45} & \frac{857}{720} \\
\end{pmatrix}
$
}

\rotatebox{90}
{\renewcommand*{\arraystretch}{1.5}
$
A_{\text{nw}} = \begin{pmatrix}
0 & 0 & 0 & 0 & 0 & 0 & 0 & 0 & 0 & 0 & 0 & 0 & 0 & 0 & 0 & 0 \\
0 & 0 & 0 & 0 & 0 & 0 & 0 & 0 & 0 & 0 & 0 & 0 & 0 & 0 & 0 & 0 \\
0 & 0 & 0 & 0 & 0 & 0 & 0 & 0 & 0 & 0 & 0 & 0 & 0 & 0 & 0 & 0 \\
0 & 0 & 0 & 0 & 0 & 0 & 0 & 0 & 0 & 0 & 0 & 0 & 0 & 0 & 0 & 0 \\
0 & 0 & 0 & 0 & \frac{1817}{720} & \frac{-229}{45} & \frac{1787}{720} & 0 & \frac{-293}{90} & \frac{1187}{180} & \frac{-571}{180} & 0 & \frac{827}{720} & \frac{-421}{180} & \frac{797}{720} & 0 \\
0 & 0 & 0 & 0 & \frac{-229}{45} & \frac{1141}{90} & \frac{-623}{90} & 0 & \frac{1187}{180} & \frac{-1547}{90} & \frac{1667}{180} & 0 & \frac{-421}{180} & \frac{278}{45} & \frac{-571}{180} & 0 \\
0 & 0 & 0 & 0 & \frac{1787}{720} & \frac{-623}{90} & \frac{3497}{720} & 0 & \frac{-571}{180} & \frac{1667}{180} & \frac{-623}{90} & 0 & \frac{797}{720} & \frac{-571}{180} & \frac{1787}{720} & 0 \\
0 & 0 & 0 & 0 & 0 & 0 & 0 & 0 & 0 & 0 & 0 & 0 & 0 & 0 & 0 & 0 \\
0 & 0 & 0 & 0 & \frac{-293}{90} & \frac{1187}{180} & \frac{-571}{180} & 0 & \frac{571}{90} & \frac{-1157}{90} & \frac{278}{45} & 0 & \frac{-109}{45} & \frac{887}{180} & \frac{-421}{180} & 0 \\
0 & 0 & 0 & 0 & \frac{1187}{180} & \frac{-1547}{90} & \frac{1667}{180} & 0 & \frac{-1157}{90} & \frac{1472}{45} & \frac{-1547}{90} & 0 & \frac{887}{180} & \frac{-1157}{90} & \frac{1187}{180} & 0 \\
0 & 0 & 0 & 0 & \frac{-571}{180} & \frac{1667}{180} & \frac{-623}{90} & 0 & \frac{278}{45} & \frac{-1547}{90} & \frac{1141}{90} & 0 & \frac{-421}{180} & \frac{1187}{180} & \frac{-229}{45} & 0 \\
0 & 0 & 0 & 0 & 0 & 0 & 0 & 0 & 0 & 0 & 0 & 0 & 0 & 0 & 0 & 0 \\
0 & 0 & 0 & 0 & \frac{827}{720} & \frac{-421}{180} & \frac{797}{720} & 0 & \frac{-109}{45} & \frac{887}{180} & \frac{-421}{180} & 0 & \frac{857}{720} & \frac{-109}{45} & \frac{827}{720} & 0 \\
0 & 0 & 0 & 0 & \frac{-421}{180} & \frac{278}{45} & \frac{-571}{180} & 0 & \frac{887}{180} & \frac{-1157}{90} & \frac{1187}{180} & 0 & \frac{-109}{45} & \frac{571}{90} & \frac{-293}{90} & 0 \\
0 & 0 & 0 & 0 & \frac{797}{720} & \frac{-571}{180} & \frac{1787}{720} & 0 & \frac{-421}{180} & \frac{1187}{180} & \frac{-229}{45} & 0 & \frac{827}{720} & \frac{-293}{90} & \frac{1817}{720} & 0 \\
0 & 0 & 0 & 0 & 0 & 0 & 0 & 0 & 0 & 0 & 0 & 0 & 0 & 0 & 0 & 0 \\
\end{pmatrix},
$
}

\rotatebox{90}
{\renewcommand*{\arraystretch}{1.5}
$
A_{\text{sw}} = \begin{pmatrix}
\frac{857}{720} & \frac{-109}{45} & \frac{827}{720} & 0 & \frac{-109}{45} & \frac{887}{180} & \frac{-421}{180} & 0 & \frac{827}{720} & \frac{-421}{180} & \frac{797}{720} & 0 & 0 & 0 & 0 & 0 \\
\frac{-109}{45} & \frac{571}{90} & \frac{-293}{90} & 0 & \frac{887}{180} & \frac{-1157}{90} & \frac{1187}{180} & 0 & \frac{-421}{180} & \frac{278}{45} & \frac{-571}{180} & 0 & 0 & 0 & 0 & 0 \\
\frac{827}{720} & \frac{-293}{90} & \frac{1817}{720} & 0 & \frac{-421}{180} & \frac{1187}{180} & \frac{-229}{45} & 0 & \frac{797}{720} & \frac{-571}{180} & \frac{1787}{720} & 0 & 0 & 0 & 0 & 0 \\
0 & 0 & 0 & 0 & 0 & 0 & 0 & 0 & 0 & 0 & 0 & 0 & 0 & 0 & 0 & 0 \\
\frac{-109}{45} & \frac{887}{180} & \frac{-421}{180} & 0 & \frac{571}{90} & \frac{-1157}{90} & \frac{278}{45} & 0 & \frac{-293}{90} & \frac{1187}{180} & \frac{-571}{180} & 0 & 0 & 0 & 0 & 0 \\
\frac{887}{180} & \frac{-1157}{90} & \frac{1187}{180} & 0 & \frac{-1157}{90} & \frac{1472}{45} & \frac{-1547}{90} & 0 & \frac{1187}{180} & \frac{-1547}{90} & \frac{1667}{180} & 0 & 0 & 0 & 0 & 0 \\
\frac{-421}{180} & \frac{1187}{180} & \frac{-229}{45} & 0 & \frac{278}{45} & \frac{-1547}{90} & \frac{1141}{90} & 0 & \frac{-571}{180} & \frac{1667}{180} & \frac{-623}{90} & 0 & 0 & 0 & 0 & 0 \\
0 & 0 & 0 & 0 & 0 & 0 & 0 & 0 & 0 & 0 & 0 & 0 & 0 & 0 & 0 & 0 \\
\frac{827}{720} & \frac{-421}{180} & \frac{797}{720} & 0 & \frac{-293}{90} & \frac{1187}{180} & \frac{-571}{180} & 0 & \frac{1817}{720} & \frac{-229}{45} & \frac{1787}{720} & 0 & 0 & 0 & 0 & 0 \\
\frac{-421}{180} & \frac{278}{45} & \frac{-571}{180} & 0 & \frac{1187}{180} & \frac{-1547}{90} & \frac{1667}{180} & 0 & \frac{-229}{45} & \frac{1141}{90} & \frac{-623}{90} & 0 & 0 & 0 & 0 & 0 \\
\frac{797}{720} & \frac{-571}{180} & \frac{1787}{720} & 0 & \frac{-571}{180} & \frac{1667}{180} & \frac{-623}{90} & 0 & \frac{1787}{720} & \frac{-623}{90} & \frac{3497}{720} & 0 & 0 & 0 & 0 & 0 \\
0 & 0 & 0 & 0 & 0 & 0 & 0 & 0 & 0 & 0 & 0 & 0 & 0 & 0 & 0 & 0 \\
0 & 0 & 0 & 0 & 0 & 0 & 0 & 0 & 0 & 0 & 0 & 0 & 0 & 0 & 0 & 0 \\
0 & 0 & 0 & 0 & 0 & 0 & 0 & 0 & 0 & 0 & 0 & 0 & 0 & 0 & 0 & 0 \\
0 & 0 & 0 & 0 & 0 & 0 & 0 & 0 & 0 & 0 & 0 & 0 & 0 & 0 & 0 & 0 \\
0 & 0 & 0 & 0 & 0 & 0 & 0 & 0 & 0 & 0 & 0 & 0 & 0 & 0 & 0 & 0 \\
\end{pmatrix},
$
}

\rotatebox{90}
{\renewcommand*{\arraystretch}{1.5}
$
A_{\text{se}} = \begin{pmatrix}
0 & 0 & 0 & 0 & 0 & 0 & 0 & 0 & 0 & 0 & 0 & 0 & 0 & 0 & 0 & 0 \\
0 & \frac{1817}{720} & \frac{-293}{90} & \frac{827}{720} & 0 & \frac{-229}{45} & \frac{1187}{180} & \frac{-421}{180} & 0 & \frac{1787}{720} & \frac{-571}{180} & \frac{797}{720} & 0 & 0 & 0 & 0 \\
0 & \frac{-293}{90} & \frac{571}{90} & \frac{-109}{45} & 0 & \frac{1187}{180} & \frac{-1157}{90} & \frac{887}{180} & 0 & \frac{-571}{180} & \frac{278}{45} & \frac{-421}{180} & 0 & 0 & 0 & 0 \\
0 & \frac{827}{720} & \frac{-109}{45} & \frac{857}{720} & 0 & \frac{-421}{180} & \frac{887}{180} & \frac{-109}{45} & 0 & \frac{797}{720} & \frac{-421}{180} & \frac{827}{720} & 0 & 0 & 0 & 0 \\
0 & 0 & 0 & 0 & 0 & 0 & 0 & 0 & 0 & 0 & 0 & 0 & 0 & 0 & 0 & 0 \\
0 & \frac{-229}{45} & \frac{1187}{180} & \frac{-421}{180} & 0 & \frac{1141}{90} & \frac{-1547}{90} & \frac{278}{45} & 0 & \frac{-623}{90} & \frac{1667}{180} & \frac{-571}{180} & 0 & 0 & 0 & 0 \\
0 & \frac{1187}{180} & \frac{-1157}{90} & \frac{887}{180} & 0 & \frac{-1547}{90} & \frac{1472}{45} & \frac{-1157}{90} & 0 & \frac{1667}{180} & \frac{-1547}{90} & \frac{1187}{180} & 0 & 0 & 0 & 0 \\
0 & \frac{-421}{180} & \frac{887}{180} & \frac{-109}{45} & 0 & \frac{278}{45} & \frac{-1157}{90} & \frac{571}{90} & 0 & \frac{-571}{180} & \frac{1187}{180} & \frac{-293}{90} & 0 & 0 & 0 & 0 \\
0 & 0 & 0 & 0 & 0 & 0 & 0 & 0 & 0 & 0 & 0 & 0 & 0 & 0 & 0 & 0 \\
0 & \frac{1787}{720} & \frac{-571}{180} & \frac{797}{720} & 0 & \frac{-623}{90} & \frac{1667}{180} & \frac{-571}{180} & 0 & \frac{3497}{720} & \frac{-623}{90} & \frac{1787}{720} & 0 & 0 & 0 & 0 \\
0 & \frac{-571}{180} & \frac{278}{45} & \frac{-421}{180} & 0 & \frac{1667}{180} & \frac{-1547}{90} & \frac{1187}{180} & 0 & \frac{-623}{90} & \frac{1141}{90} & \frac{-229}{45} & 0 & 0 & 0 & 0 \\
0 & \frac{797}{720} & \frac{-421}{180} & \frac{827}{720} & 0 & \frac{-571}{180} & \frac{1187}{180} & \frac{-293}{90} & 0 & \frac{1787}{720} & \frac{-229}{45} & \frac{1817}{720} & 0 & 0 & 0 & 0 \\
0 & 0 & 0 & 0 & 0 & 0 & 0 & 0 & 0 & 0 & 0 & 0 & 0 & 0 & 0 & 0 \\
0 & 0 & 0 & 0 & 0 & 0 & 0 & 0 & 0 & 0 & 0 & 0 & 0 & 0 & 0 & 0 \\
0 & 0 & 0 & 0 & 0 & 0 & 0 & 0 & 0 & 0 & 0 & 0 & 0 & 0 & 0 & 0 \\
0 & 0 & 0 & 0 & 0 & 0 & 0 & 0 & 0 & 0 & 0 & 0 & 0 & 0 & 0 & 0 \\
\end{pmatrix}.
$
}

The matrix $A_{\text{opt}}$ appearing in \eqref{eq:ind2dDati} is splitted into two matrices $A_{\text{opt}}^{(1)}$ and $A_{\text{opt}}^{(2)}$ reporting the first eight columns and the second eight columns of $A_{\text{opt}}$.

\rotatebox{90}
{\renewcommand*{\arraystretch}{1.5}
$
A_{\text{opt}}^{(1)} = \begin{pmatrix}
\frac{2903}{1575} & \frac{-182131}{37800} & \frac{10859}{2700} & \frac{-13889}{12600} & \frac{-182131}{37800} & \frac{60463}{4800} & \frac{-58881}{5600} & \frac{871553}{302400} \\
\frac{-182131}{37800} & \frac{48799}{3150} & \frac{-58687}{4200} & \frac{10859}{2700} & \frac{60463}{4800} & \frac{-340087}{8400} & \frac{1226693}{33600} & \frac{-58881}{5600} \\
\frac{10859}{2700} & \frac{-58687}{4200} & \frac{48799}{3150} & \frac{-182131}{37800} & \frac{-58881}{5600} & \frac{1226693}{33600} & \frac{-340087}{8400} & \frac{60463}{4800} \\
\frac{-13889}{12600} & \frac{10859}{2700} & \frac{-182131}{37800} & \frac{2903}{1575} & \frac{871553}{302400} & \frac{-58881}{5600} & \frac{60463}{4800} & \frac{-182131}{37800} \\
\frac{-182131}{37800} & \frac{60463}{4800} & \frac{-58881}{5600} & \frac{871553}{302400} & \frac{48799}{3150} & \frac{-340087}{8400} & \frac{426623}{12600} & \frac{-11137}{1200} \\
\frac{60463}{4800} & \frac{-340087}{8400} & \frac{1226693}{33600} & \frac{-58881}{5600} & \frac{-340087}{8400} & \frac{815527}{6300} & \frac{-2958593}{25200} & \frac{426623}{12600} \\
\frac{-58881}{5600} & \frac{1226693}{33600} & \frac{-340087}{8400} & \frac{60463}{4800} & \frac{426623}{12600} & \frac{-2958593}{25200} & \frac{815527}{6300} & \frac{-340087}{8400} \\
\frac{871553}{302400} & \frac{-58881}{5600} & \frac{60463}{4800} & \frac{-182131}{37800} & \frac{-11137}{1200} & \frac{426623}{12600} & \frac{-340087}{8400} & \frac{48799}{3150} \\
\frac{10859}{2700} & \frac{-58881}{5600} & \frac{221047}{25200} & \frac{-363499}{151200} & \frac{-58687}{4200} & \frac{1226693}{33600} & \frac{-1536757}{50400} & \frac{841823}{100800} \\
\frac{-58881}{5600} & \frac{426623}{12600} & \frac{-1536757}{50400} & \frac{221047}{25200} & \frac{1226693}{33600} & \frac{-2958593}{25200} & \frac{10709107}{100800} & \frac{-1536757}{50400} \\
\frac{221047}{25200} & \frac{-1536757}{50400} & \frac{426623}{12600} & \frac{-58881}{5600} & \frac{-1536757}{50400} & \frac{10709107}{100800} & \frac{-2958593}{25200} & \frac{1226693}{33600} \\
\frac{-363499}{151200} & \frac{221047}{25200} & \frac{-58881}{5600} & \frac{10859}{2700} & \frac{841823}{100800} & \frac{-1536757}{50400} & \frac{1226693}{33600} & \frac{-58687}{4200} \\
\frac{-13889}{12600} & \frac{871553}{302400} & \frac{-363499}{151200} & \frac{66427}{100800} & \frac{10859}{2700} & \frac{-58881}{5600} & \frac{221047}{25200} & \frac{-363499}{151200} \\
\frac{871553}{302400} & \frac{-11137}{1200} & \frac{841823}{100800} & \frac{-363499}{151200} & \frac{-58881}{5600} & \frac{426623}{12600} & \frac{-1536757}{50400} & \frac{221047}{25200} \\
\frac{-363499}{151200} & \frac{841823}{100800} & \frac{-11137}{1200} & \frac{871553}{302400} & \frac{221047}{25200} & \frac{-1536757}{50400} & \frac{426623}{12600} & \frac{-58881}{5600} \\
\frac{66427}{100800} & \frac{-363499}{151200} & \frac{871553}{302400} & \frac{-13889}{12600} & \frac{-363499}{151200} & \frac{221047}{25200} & \frac{-58881}{5600} & \frac{10859}{2700} \\
\end{pmatrix},
$
}

\rotatebox{90}
{\renewcommand*{\arraystretch}{1.5}
$
A_{\text{opt}}^{(2)} = \begin{pmatrix}
\frac{10859}{2700} & \frac{-58881}{5600} & \frac{221047}{25200} & \frac{-363499}{151200} & \frac{-13889}{12600} & \frac{871553}{302400} & \frac{-363499}{151200} & \frac{66427}{100800} \\
\frac{-58881}{5600} & \frac{426623}{12600} & \frac{-1536757}{50400} & \frac{221047}{25200} & \frac{871553}{302400} & \frac{-11137}{1200} & \frac{841823}{100800} & \frac{-363499}{151200} \\
\frac{221047}{25200} & \frac{-1536757}{50400} & \frac{426623}{12600} & \frac{-58881}{5600} & \frac{-363499}{151200} & \frac{841823}{100800} & \frac{-11137}{1200} & \frac{871553}{302400} \\
\frac{-363499}{151200} & \frac{221047}{25200} & \frac{-58881}{5600} & \frac{10859}{2700} & \frac{66427}{100800} & \frac{-363499}{151200} & \frac{871553}{302400} & \frac{-13889}{12600} \\
\frac{-58687}{4200} & \frac{1226693}{33600} & \frac{-1536757}{50400} & \frac{841823}{100800} & \frac{10859}{2700} & \frac{-58881}{5600} & \frac{221047}{25200} & \frac{-363499}{151200} \\
\frac{1226693}{33600} & \frac{-2958593}{25200} & \frac{10709107}{100800} & \frac{-1536757}{50400} & \frac{-58881}{5600} & \frac{426623}{12600} & \frac{-1536757}{50400} & \frac{221047}{25200} \\
\frac{-1536757}{50400} & \frac{10709107}{100800} & \frac{-2958593}{25200} & \frac{1226693}{33600} & \frac{221047}{25200} & \frac{-1536757}{50400} & \frac{426623}{12600} & \frac{-58881}{5600} \\
\frac{841823}{100800} & \frac{-1536757}{50400} & \frac{1226693}{33600} & \frac{-58687}{4200} & \frac{-363499}{151200} & \frac{221047}{25200} & \frac{-58881}{5600} & \frac{10859}{2700} \\
\frac{48799}{3150} & \frac{-340087}{8400} & \frac{426623}{12600} & \frac{-11137}{1200} & \frac{-182131}{37800} & \frac{60463}{4800} & \frac{-58881}{5600} & \frac{871553}{302400} \\
\frac{-340087}{8400} & \frac{815527}{6300} & \frac{-2958593}{25200} & \frac{426623}{12600} & \frac{60463}{4800} & \frac{-340087}{8400} & \frac{1226693}{33600} & \frac{-58881}{5600} \\
\frac{426623}{12600} & \frac{-2958593}{25200} & \frac{815527}{6300} & \frac{-340087}{8400} & \frac{-58881}{5600} & \frac{1226693}{33600} & \frac{-340087}{8400} & \frac{60463}{4800} \\
\frac{-11137}{1200} & \frac{426623}{12600} & \frac{-340087}{8400} & \frac{48799}{3150} & \frac{871553}{302400} & \frac{-58881}{5600} & \frac{60463}{4800} & \frac{-182131}{37800} \\
\frac{-182131}{37800} & \frac{60463}{4800} & \frac{-58881}{5600} & \frac{871553}{302400} & \frac{2903}{1575} & \frac{-182131}{37800} & \frac{10859}{2700} & \frac{-13889}{12600} \\
\frac{60463}{4800} & \frac{-340087}{8400} & \frac{1226693}{33600} & \frac{-58881}{5600} & \frac{-182131}{37800} & \frac{48799}{3150} & \frac{-58687}{4200} & \frac{10859}{2700} \\
\frac{-58881}{5600} & \frac{1226693}{33600} & \frac{-340087}{8400} & \frac{60463}{4800} & \frac{10859}{2700} & \frac{-58687}{4200} & \frac{48799}{3150} & \frac{-182131}{37800} \\
\frac{871553}{302400} & \frac{-58881}{5600} & \frac{60463}{4800} & \frac{-182131}{37800} & \frac{-13889}{12600} & \frac{10859}{2700} & \frac{-182131}{37800} & \frac{2903}{1575} \\
\end{pmatrix}.
$
}

\end{document}